\documentclass{amsart}
\usepackage[T1]{fontenc}
\usepackage{lmodern}
\usepackage{amsmath,amssymb,amsthm}
\usepackage[hidelinks]{hyperref}
\usepackage[margin=2.6cm]{geometry}

\usepackage{tikz}
\tikzset{dot/.style={circle,fill,inner sep=1pt},wdot/.style={circle,draw,fill=white,inner sep=1pt},lab/.style={font=\scriptsize,inner sep=1.5pt}}

\newtheorem{theorem}{Theorem}[section]
\newtheorem{lemma}[theorem]{Lemma}
\newtheorem{proposition}[theorem]{Proposition}
\newtheorem{corollary}[theorem]{Corollary}
\newtheorem{conjecture}[theorem]{Conjecture}
\newtheorem*{theoremA}{Theorem A}
\newtheorem*{theoremB}{Theorem B}
\newtheorem*{corollaryC}{Corollary C}
\theoremstyle{definition}
\newtheorem{definition}[theorem]{Definition}
\newtheorem{example}[theorem]{Example}

\DeclareMathOperator{\rk}{rank}
\DeclareMathOperator{\htop}{ht}
\DeclareMathOperator{\Min}{Min}
\DeclareMathOperator{\col}{col}
\newcommand{\cC}{\mathcal C}
\newcommand{\cD}{\mathcal D}
\newcommand{\cE}{\mathcal E}
\newcommand{\cN}{\mathcal N}
\newcommand{\Q}{\mathbb Q}
\newcommand{\Z}{\mathbb Z}
\newcommand{\N}{\mathbb N}
\newcommand{\R}{\mathbb R}
\newcommand{\one}{\mathbf 1}
\newcommand{\dm}{\diamond}
\newcommand{\lt}{\underline{t}}
\newcommand{\ut}{\overline{t}}

\title[Join-meet ideals of modular lattices with one diamond]{Minimal primes and radicality of join-meet ideals of modular lattices with one diamond}

\author[T.~Hibi]{Takayuki Hibi}
\address[Takayuki Hibi]
{Department of Pure and Applied Mathematics, 
	Graduate School of Information Science and Technology, 
	Osaka University, 
	Suita, Osaka 565-0871, Japan}
\email{hibi@math.sci.osaka-u.ac.jp}

\author[F. Navarra]{Francesco Navarra}      
\address[Francesco Navarra]{Sabanci University, Faculty of Engineering and Natural Sciences, Orta Mahalle, Tuzla 34956, Istanbul, Turkey}	
\email{francesco.navarra@sabanciuniv.edu}



\thanks{The second author is supported by Scientific and Technological Research Council of Turkey T\"UB\.{I}TAK under the Grants 126F175; he is thankful to T\"UB\.{I}TAK for their support. He is also a member of GNSAGA Indam and he acknowledges their support.}

\begin{document}

\begin{abstract}
Let $L$ be a finite modular lattice with exactly one diamond, that is, exactly one interval of length two with at least three atoms. We describe the minimal primes of the join-meet ideal $I_L$ in terms of maximal closed distributive sublattices of $L$, and we show that $I_L$ has height $|L|-\rk L$. Moreover, we characterize the radicality of $I_L$, over any field, in terms of the complements of the atoms of the diamond in their fibres. Finally, we prove that $I_L$ is unmixed if and only if all its minimal primes have the same height, if and only if the rank of every maximal closed distributive sublattice through an atom of the diamond is $\rk L-1$; in this case $I_L$ is radical.
\end{abstract}

\maketitle

\section{Introduction}

Let $L$ be a finite lattice and $K$ a field. The join-meet ideal of $L$ is
\[
I_L=(x_ax_b-x_{a\wedge b}x_{a\vee b}:a,b\in L\text{ incomparable})\subseteq S_L=K[x_a:a\in L].
\]
Join-meet ideals have been widely studied; see \cite[Chapter 6]{HHO} and the references therein. If $L$ is distributive, that is, if $L$ has no sublattice isomorphic to the pentagon $N_5$ or to $M_3$, much is known: $I_L$ is the toric ideal of the Hibi ring of $L$, which is a normal Cohen--Macaulay domain \cite{H}. If $L$ is modular but not distributive, that is, if $L$ has a sublattice isomorphic to $M_3$ but none isomorphic to $N_5$, most algebraic properties of $I_L$ are still unknown. The first systematic study of this case is due to Ene and Hibi \cite{EH}: they proved that no initial ideal of $I_L$ is squarefree as conjectured in \cite{HH12}, found a class of modular non-distributive lattices whose join-meet ideal is radical, and described the minimal primes of radical join-meet ideals in terms of admissible sets. For modular non-distributive lattices, the questions of when $I_L$ is radical or unmixed are still open, and they seem to be difficult in general. Inspired by \cite{EH}, in this paper we answer these questions, and we determine the minimal primes of $I_L$, when $L$ is a modular lattice with exactly one diamond; since every modular non-distributive lattice has at least one diamond, this is the first case beyond the distributive one.

Let $L$ be a finite modular lattice with exactly one diamond, that is, exactly one interval $[p,q]$ of length two with at least three atoms. Let $A_\dm$ be the set of atoms of $[p,q]$, $n=|L|$ and $r=\rk L$, where the rank of a lattice, or of a sublattice, is the largest length of a chain contained in it. A nonempty set $W\subseteq L$ is \emph{closed} if, for all incomparable $a,b\in L$, we have $a,b\in W$ if and only if $a\wedge b,a\vee b\in W$; equivalently, $L\setminus W$ is an admissible set in the sense of \cite[Definition 2.2]{EH}. Closed sets are sublattices of $L$. A \emph{closed distributive sublattice} of $L$ is a closed set which is a distributive lattice, and it is \emph{maximal} if no other closed distributive sublattice of $L$ contains it. We denote by $\cD(L)$ the set of maximal closed distributive sublattices of $L$ containing an atom of $[p,q]$, and by $\cN(L)$ the set of maximal closed distributive sublattices $W$ of $L$ with $W\cap[p,q]=\emptyset$ and $\rk W=r-1$; see Definition~\ref{def:DN} and Example~\ref{ex:DN}. For a closed set $W$, let $P_W$ be the ideal generated by the variables $x_a$ with $a\notin W$ and by the saturation $I_W:(\prod_{a\in W}x_a)^\infty$ of the join-meet ideal of the lattice $W$; it is a prime ideal containing $I_L$ (Proposition~\ref{prop:minprimes}). Our first result is a complete description of the minimal primes of $I_L$.

\begin{theoremA}
	Let $L$ be a finite modular lattice with exactly one diamond $[p,q]$. Then
	\[
	\Min(I_L)=\{P_L\}\cup\{P_D : D\in\cD(L)\}\cup\{P_W : W\in\cN(L)\},
	\]
	and distinct sets give distinct primes. Moreover $\htop P_L=\htop P_W=n-r$ for $W\in\cN(L)$, and $\htop P_D=n-1-\rk D\ge n-r$ for $D\in\cD(L)$. In particular $\htop I_L=n-r$.
\end{theoremA}

For $c\in A_\dm$, the elements $x$ with $(x\vee p)\wedge q=c$ form an interval $[m_c,M_c]$, the \emph{fibre} of $c$ (Lemma~\ref{lem:rho}), and we say that $c$ is \emph{complemented in its fibre} if there is $c'\in[m_c,M_c]$ with $c\wedge c'=m_c$ and $c\vee c'=M_c$. Intuitively, the map $x\mapsto(x\vee p)\wedge q$ projects $L$ onto its diamond, and the fibre of $c$ is the part of $L$ lying over $c$.\\
It turns out that the radicality of $I_L$ depends only on whether the atoms of the diamond are complemented in their fibres, and this leads to a characterization of unmixed join-meet ideals.

\begin{theoremB}
	Let $L$ be a finite modular lattice with exactly one diamond. Then $I_L$ is radical if and only if every atom of the diamond is complemented in its fibre.
\end{theoremB}

\begin{corollaryC}
	Let $L$ be a finite modular lattice with exactly one diamond. The following are equivalent:
	\begin{enumerate}
		\item[(a)] $I_L$ is unmixed;
		\item[(b)] all minimal primes of $I_L$ have the same height;
		\item[(c)] $\rk D=r-1$ for every $D\in\cD(L)$.
	\end{enumerate}
	If they hold, then $I_L$ is radical.
\end{corollaryC}

The two theorems are linked by Proposition~\ref{prop:local}: an atom that is not complemented in its fibre forces some $D\in\cD(L)$ of rank at most $r-2$, while the converse fails. Figure~\ref{fig:examples} and Examples~\ref{ex:A}--\ref{ex:L12} illustrate Theorem B and Corollary C. Since no initial ideal of $I_L$ is squarefree when $L$ is modular and not distributive \cite[Theorem 1.3]{EH}, Gr\"obner degenerations cannot prove Theorem B, and our proof is combinatorial: $I_L$ is radical if and only if, in every multidegree, at most one class of monomials under the moves $x_ax_b\leftrightarrow x_{a\wedge b}x_{a\vee b}$ is not supported on a closed distributive sublattice (Theorem~\ref{thm:radcrit}).

The paper is organized as follows. Section~\ref{sec:prelim} recalls the background. In Section~\ref{sec:toric}, for an arbitrary finite modular lattice, we identify the saturation of $I_L$ with a toric ideal, we show that the minimal primes of $I_L$ are among the primes $P_W$, and we give a criterion for inclusions among them. Section~\ref{sec:structure} describes the structure of a modular lattice with one diamond; we use it in Section~\ref{sec:inclusions} to compare the primes $P_W$, and in Section~\ref{sec:main} to prove Theorem A. Section~\ref{sec:radical} is devoted to Theorem B and Corollary C, and Section~\ref{sec:final} to examples and conjectures. In particular, Examples~\ref{ex:A} and~\ref{ex:Z} show that unmixedness cannot be characterized by forbidden sublattices: the lattice $A$ is an interval of the lattice $Z$, and $I_Z$ is unmixed while $I_A$ is not.

\section{Preliminaries}\label{sec:prelim}

\subsection{Lattices}
All posets and lattices are finite. A \emph{chain} of a poset $\Pi$ is a totally ordered subset of $\Pi$, its \emph{length} is its cardinality minus one, the \emph{rank} $\rk\Pi$ of $\Pi$ is the largest length of a chain of $\Pi$. A subset $U\subseteq\Pi$ is an \emph{order ideal} if $y\in U$ whenever $y\le x$ and $x\in U$. We denote by $\mathcal J(\Pi)$ the set of order ideals of $\Pi$, ordered by inclusion; it is a sublattice of the Boolean lattice $2^\Pi$ of all subsets of $\Pi$.

A \emph{lattice} is a poset $L$ in which any two elements $a,b$ have a greatest lower bound $a\wedge b$, their \emph{meet}, and a least upper bound $a\vee b$, their \emph{join}; we denote by $\hat0$ and $\hat1$ the least and the greatest element of $L$, and $\rk L$ is the largest length of a chain of $L$. A \emph{sublattice} of $L$ is a nonempty subset closed under meet and join; it is a lattice with the order of $L$, and its rank is defined as above. For $a\le b$ we put $[a,b]=\{x\in L:a\le x\le b\}$, ${\downarrow}a=[\hat0,a]$ and ${\uparrow}a=[a,\hat1]$; these \emph{intervals} are sublattices, and every nonempty convex sublattice is an interval. If $a<b$ and $[a,b]=\{a,b\}$, then $b$ \emph{covers} $a$, we write $a\lessdot b$, and $[a,b]$ is a \emph{prime interval}. The \emph{length} of $[a,b]$ is its rank, and the \emph{atoms} of $[a,b]$ are the elements of $[a,b]$ covering $a$. An element $c'\in[a,b]$ is a \emph{complement} of $c\in[a,b]$ in $[a,b]$ if $c\wedge c'=a$ and $c\vee c'=b$. The lattice $L$ is \emph{graded} if all its maximal chains have the same length, and then $\rk a$ is the length of ${\downarrow}a$. An element $a\ne\hat0$ is \emph{join-irreducible} if $a=b\vee c$ implies $a\in\{b,c\}$. The dual lattice $L^{\mathrm{op}}$ has the reverse order.

A lattice is \emph{modular} if $a\vee(x\wedge b)=(a\vee x)\wedge b$ whenever $a\le b$, and \emph{distributive} if $x\wedge(y\vee z)=(x\wedge y)\vee(x\wedge z)$ for all $x,y,z$. Both properties pass to sublattices, products, homomorphic images and dual lattices. By Birkhoff's theorem \cite[Theorem~9.1.7]{HH11}, if $L$ is distributive and $\Pi$ is the poset of its join-irreducible elements, then $a\mapsto\{b\in\Pi:b\le a\}$ is an isomorphism of lattices from $L$ onto $\mathcal J(\Pi)$; moreover, $\rk\mathcal J(\Pi)=|\Pi|$. A lattice is modular if and only if it has no sublattice isomorphic to the pentagon $N_5$ (Dedekind \cite{De}, see also \cite[pp.~65--66]{B}), and a modular lattice is distributive if and only if it has no sublattice isomorphic to the diamond $M_3$ (Birkhoff \cite[p.~152]{B}); see also \cite[Theorem 6.10]{HHO}. We call \emph{diamond} of $L$ an interval of length two with at least three atoms; a modular lattice with a diamond is not distributive.


\subsection{Projectivity and the Jordan--H\"older valuation}
Let $L$ be modular. For prime intervals we write $[a,b]\nearrow[c,d]$, or $[c,d]\searrow[a,b]$, if $b\wedge c=a$ and $b\vee c=d$, and we call $[a,b]$ and $[c,d]$ \emph{projective}, $[a,b]\approx[c,d]$, if they are connected by a sequence of such transpositions. These are the perspectivities and the projectivity of \cite[Section I.3.5]{G}, restricted to prime intervals; this is harmless, since perspective intervals of a modular lattice are isomorphic \cite[Theorem 348]{G}. Let $\cC(L)$ be the set of classes of $\approx$ and $\kappa(L)=|\cC(L)|$. For instance, the six prime intervals of $M_3$ form a single class.

Let $x\le y$. By the Jordan--H\"older theorem \cite[Theorems 374 and 377]{G}, which applies because modular lattices are semimodular \cite[Theorem 350]{G}, any two maximal chains of $[x,y]$ have the same length, and their prime intervals correspond bijectively under $\approx$. So $L$ is graded, and for $C\in\cC(L)$ the number $m_C(x,y)$ of prime intervals of class $C$ on a maximal chain of $[x,y]$ does not depend on the chain. We put $m(x,y)=(m_C(x,y))_C\in\N^{\cC(L)}$ and call
\[
v:L\to\N^{\cC(L)},\qquad v(x)=m(\hat0,x),
\]
the \emph{Jordan--H\"older valuation} of $L$. Thus $v(y)-v(x)=m(x,y)$ for $x\le y$, the entries of $m(x,y)$ add up to the length of $[x,y]$, and $v(b)=v(a)+e_C$ if $[a,b]\in C$, where $e_C$ is a unit vector. The following facts are consequences of the isomorphism theorem for modular lattices \cite[Theorem 348]{G} and of the Jordan--H\"older theorem.

\begin{lemma}\label{lem:valuation}
	Let $L$ be modular and $a,b\in L$. Then $m(a\wedge b,a)=m(b,a\vee b)$ and $v(a)+v(b)=v(a\wedge b)+v(a\vee b)$. Every class $C$ meets every maximal chain of $L$, and $\sum_Cm_C(\hat0,\hat1)=\rk L$.
\end{lemma}

In particular $\rk a+\rk b=\rk(a\wedge b)+\rk(a\vee b)$ for all $a,b\in L$.

\begin{lemma}\label{lem:kappa}
	Let $T$ be a modular lattice. Then $\kappa(T)\le\rk T$, with equality if and only if $T$ is distributive. If $T$ is distributive, every class of $T$ contains exactly one prime interval of each maximal chain.
\end{lemma}

\begin{proof}
	By Lemma~\ref{lem:valuation}, $\kappa(T)\le\sum_Cm_C(\hat0,\hat1)=\rk T$, with equality if and only if every class occurs once on each maximal chain. If $T$ is not distributive, it has a sublattice $\{o,x,y,z,i\}\cong M_3$, and Lemma~\ref{lem:valuation} gives $m(o,x)=m(y,i)=m(o,z)=m(x,i)$, so a maximal chain through $o<x<i$ meets some class twice. If $T$ is distributive, we may assume $T=\mathcal J(\Pi)$; then a cover adds one element of $\Pi$, and if $[a,b]\nearrow[c,d]$, then $d\setminus c=(b\cup c)\setminus c=b\setminus a$; so projective covers add the same element, and since a maximal chain adds each element of $\Pi$ once, it meets each class once.
\end{proof}

\subsection{Join-meet ideals and closed sets}\label{subsec:closed}
Let $K$ be a field. The \emph{join-meet ideal} of a finite lattice $L$ is the ideal $I_L\subseteq S_L=K[x_a:a\in L]$ generated by the binomials $x_ax_b-x_{a\wedge b}x_{a\vee b}$ with $a,b\in L$ incomparable. If $L$ is distributive, $I_L$ is prime by Hibi's theorem \cite{H}; see also \cite[Theorem 6.21]{HHO}. We now introduce the notions needed to describe the minimal primes of $I_L$ and to study its unmixedness. Recall that an ideal is \emph{unmixed} if all its associated primes have the same height, and $\Min(I)$ denotes the set of minimal primes of an ideal $I$. Following \cite[Definition 2.2]{EH}, a subset $A\subseteq L$ is \emph{admissible} if, for every generator $x_ax_b-x_{a\wedge b}x_{a\vee b}$ of $I_L$, either both monomials or neither of them is divisible by a variable $x_c$ with $c\in A$. A nonempty set $W\subseteq L$ is \emph{closed} if $L\setminus W$ is admissible. 

\begin{figure}[htbp]
	\centering
	\begin{tikzpicture}[x=0.6cm,y=0.55cm]
		\foreach \s/\lbl/\blk/\wht in {0/{(a) $\{0,1,3,7,8\}$, not closed}/{0/below,1/left,3/right,7/right,8/above}/{2/right,4/right,5/left,6/left},5.2/{(b) $D=\{0,1,2,5,6\}\in\cD(L)$}/{0/below,1/left,2/right,5/left,6/left}/{3/right,4/right,7/right,8/above},10.4/{(c) $W=\{0,1,6,8\}\in\cN(L)$}/{0/below,1/left,6/left,8/above}/{2/right,3/right,4/right,5/left,7/right}}{
			\begin{scope}[xshift=\s cm]
				\foreach \n/\x/\y in {0/0/0,1/-1/1,2/1/1,3/1/2,4/2/2,5/0/2,6/-1/3,7/1/3,8/0/4} \coordinate (v\n) at (\x,\y);
				\foreach \a/\b in {0/1,0/2,1/5,2/5,2/3,2/4,5/6,5/7,3/7,4/7,6/8,7/8} \draw (v\a)--(v\b);
				\foreach \n/\p in \blk \node[dot,label={[lab]\p:$\n$}] at (v\n) {};
				\foreach \n/\p in \wht \node[wdot,label={[lab]\p:$\n$}] at (v\n) {};
				\node[lab] at (0.5,-1.4) {\lbl};
		\end{scope}}
	\end{tikzpicture}
	\caption{The lattice $L$ of Examples~\ref{ex:closed} and~\ref{ex:DN}, whose only diamond is $[2,7]$. In each panel the black elements form the set under consideration and the white ones its complement.}\label{fig:DN}
\end{figure}
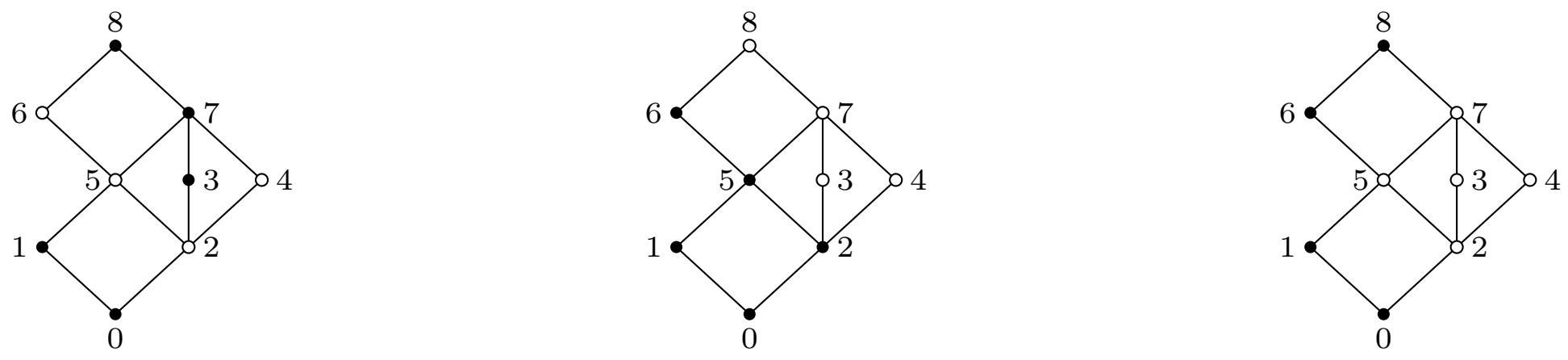

\begin{example}\label{ex:closed}
	Let $L$ be the lattice of Figure~\ref{fig:DN}, whose only diamond is $[2,7]$. The set $D=\{0,1,2,5,6\}$ of Figure~\ref{fig:DN}(b) is closed, that is, $L\setminus D=\{3,4,7,8\}$ is admissible: one checks that every generator of $I_L$ has either both monomials or neither of them divisible by one of $x_3,x_4,x_7,x_8$. For instance, in $x_6x_7-x_5x_8$ the first monomial is divisible by $x_7$ and the second by $x_8$, while in $x_1x_2-x_0x_5$ neither monomial is. On the other hand, the set $\{0,1,3,7,8\}$ of Figure~\ref{fig:DN}(a) is not closed, since its complement $\{2,4,5,6\}$ is not admissible: in the generator $x_1x_4-x_0x_7$ the first monomial is divisible by $x_4$, while the second is divisible by none of $x_2,x_4,x_5,x_6$. This set is nevertheless a sublattice of $L$; we will see in Proposition~\ref{prop:closed} that, conversely, every closed set is a sublattice.
\end{example}

The following elementary properties of closed sets will be used throughout the paper.

\begin{proposition}\label{prop:closed}
	Let $W\subseteq L$ be a nonempty set.
	\begin{enumerate}
		\item[(a)] $W$ is closed if and only if, for all incomparable $a,b\in L$, we have $a,b\in W$ if and only if $a\wedge b,a\vee b\in W$.
		\item[(b)] If $W$ is closed, then $W$ is a sublattice of $L$, and for $x\le y$ the set $W\cap[x,y]$ is either empty or a closed subset of the lattice $[x,y]$.
		\item[(c)] Every interval of $L$ is closed, and every nonempty intersection of closed sets is closed.
	\end{enumerate}
\end{proposition}

\begin{proof}
	(a) Let $A=L\setminus W$ and let $a,b\in L$ be incomparable. The monomial $x_ax_b$ is divisible by a variable $x_c$ with $c\in A$ if and only if $\{a,b\}\not\subseteq W$, and $x_{a\wedge b}x_{a\vee b}$ is divisible by such a variable if and only if $\{a\wedge b,a\vee b\}\not\subseteq W$. So $A$ is admissible exactly when the condition in (a) holds.
	
	(b) If $a,b\in W$ are comparable, then $\{a\wedge b,a\vee b\}=\{a,b\}\subseteq W$; if they are incomparable, then $a\wedge b,a\vee b\in W$ by (a). Meets, joins and comparability of elements of $[x,y]$ are the same in $[x,y]$ and in $L$, so for incomparable $a,b\in[x,y]$ the condition in (a) reads $a,b\in W\cap[x,y]$ if and only if $a\wedge b,a\vee b\in W\cap[x,y]$.
	
	(c) If $a,b$ are incomparable, then $a\wedge b<a,b<a\vee b$, so $a,b\in[x,y]$ if and only if $a\wedge b,a\vee b\in[x,y]$. The condition in (a) is clearly preserved by intersections.
\end{proof}
For a closed set $W$, let $R_W=K[x_a:a\in W]$, let $I_W\subseteq R_W$ be the join-meet ideal of the lattice $W$, and put
\[
Q_W=I_W:\Bigl(\prod_{a\in W}x_a\Bigr)^\infty\subseteq R_W\qquad\text{and}\qquad P_W=(x_a:a\in L\setminus W)+Q_WS_L.
\]
We will see in Proposition~\ref{prop:minprimes} that $P_W$ is a prime ideal containing $I_L$, and that every minimal prime of $I_L$ is of this form. The sets $W$ that actually occur are described by the following notions.

\begin{definition}\label{def:DN}
	Let $L$ be a finite modular lattice. A \emph{closed distributive sublattice} of $L$ is a closed set $D\subseteq L$ which is a distributive lattice, and it is \emph{maximal} if no other closed distributive sublattice of $L$ contains it. We denote by
	\begin{itemize}
		\item $\cD(L)$ the set of maximal closed distributive sublattices of $L$ containing an atom of some diamond of $L$;
		\item $\cN(L)$, when $L$ has exactly one diamond $[p,q]$, the set of maximal closed distributive sublattices $W$ of $L$ with $W\cap[p,q]=\emptyset$ and $\rk W=\rk L-1$.
	\end{itemize}
\end{definition}

\begin{example}\label{ex:DN}
	The lattice $L$ of Figure~\ref{fig:DN} has rank $4$. One checks that $\cN(L)$ consists only of the chain $W=\{0,1,6,8\}$ of Figure~\ref{fig:DN}(c), and that $\cD(L)$ consists of the six chains $\{0,2,3\}$, $\{0,2,4\}$, $\{0,3,8\}$, $\{0,4,8\}$, $\{3,7,8\}$, $\{4,7,8\}$, of rank $2$, and of the sublattices $\{0,1,2,5,6\}$ and $\{1,5,6,7,8\}$, of rank $3$; the first of them is shown in Figure~\ref{fig:DN}(b). By Theorem A, $I_L$ has ten minimal primes, four of height $5$ and six of height $6$.
\end{example}

\section{The toric component and the minimal primes}\label{sec:toric}

In this section $L$ is a finite modular lattice, $n=|L|$ and $r=\rk L$. We first identify the saturation $Q_L$ of $I_L$ with a toric ideal, and for this we recall some facts on lattice ideals; here a \emph{lattice} is a subgroup of $\Z^k$, not a poset. For a subgroup $\Gamma\subseteq\Z^k$, the \emph{lattice ideal} $I_\Gamma\subseteq K[x_1,\dots,x_k]$ is generated by the binomials $x^{z^+}-x^{z^-}$, $z\in\Gamma$, where $z^+$ and $z^-$ are the positive and negative parts of $z$. Such a binomial lies in $I_\Gamma$ if and only if $z\in\Gamma$, and $I_\Gamma:(x_1\cdots x_k)^\infty=I_\Gamma$ \cite[Theorem 3.20 and its proof]{HHO}; moreover, $I_\Gamma$ is prime if and only if $\Gamma$ is saturated, that is, $\Z^k/\Gamma$ is torsion-free \cite[Theorem 3.17]{HHO}. For $h_1,\dots,h_k\in\Z^d$, the kernel of the monomial map $x_i\mapsto t^{h_i}$ is the lattice ideal of the kernel of $\Z^k\to\Z^d$, $e_i\mapsto h_i$ \cite[Proposition 3.16]{HHO}, and its height is $k$ minus the rank of $h_1,\dots,h_k$ \cite[Proposition 3.1]{HHO}.

\begin{theorem}\label{thm:toric}
	Let $L$ be a finite modular lattice and $K$ any field. Then the subgroup $\Lambda_L\subseteq\Z^L$ generated by the vectors $e_a+e_b-e_{a\wedge b}-e_{a\vee b}$ is the kernel of $\Phi_L:\Z^L\to\Z\oplus\Z^{\cC(L)}$, $e_a\mapsto(1,v(a))$, and it is saturated of rank $|L|-1-\kappa(L)$. Moreover, $Q_L$ is the kernel of $\phi_L:S_L\to K[t,s_C : C\in\cC(L)]$, $x_a\mapsto t\,s^{v(a)}$; in particular $Q_L$ is prime and $\htop Q_L=|L|-1-\kappa(L)$.
\end{theorem}

\begin{proof}
	Fix a maximal chain $\hat0=c_0\lessdot\cdots\lessdot c_r=\hat1$. By Lemma~\ref{lem:valuation} every class $C$ occurs on it, and we choose $i_C$ with $[c_{i_C-1},c_{i_C}]\in C$. Put $h_C=e_{c_{i_C}}-e_{c_{i_C-1}}$, $F=\{\hat0\}\cup\{c_{i_C}:C\in\cC(L)\}$, and $b_a=e_a-e_{\hat0}-\sum_Cv_C(a)h_C$ for $a\in L\setminus F$. The $|L|$ vectors $e_{\hat0}$, $h_C$ and $b_a$ generate $\Z^L$: this is clear for $e_a$ with $a\notin F$, and for $e_{c_j}$ it follows by induction on $j$, since $e_{c_j}=e_{c_{j-1}}+h_C$ if $j=i_C$. Hence they form a basis of $\Z^L$. If $[u,u']\nearrow[w,w']$, then $(e_{u'}-e_u)-(e_{w'}-e_w)$ is a generator of $\Lambda_L$. Since projectivity is generated by transpositions, $e_{u'}-e_u\equiv h_C$ modulo $\Lambda_L$ for every prime interval $[u,u']\in C$, and summing along a maximal chain of ${\downarrow}a$ gives $b_a\in\Lambda_L$. Now $\Phi_L(e_{\hat0})=(1,0)$, $\Phi_L(h_C)=(0,e_C)$ and $\Phi_L(b_a)=0$, and $\Lambda_L\subseteq\ker\Phi_L$ by Lemma~\ref{lem:valuation}; therefore $\Lambda_L=\ker\Phi_L$ is the direct summand of $\Z^L$ spanned by the vectors $b_a$. Finally, $I_L\subseteq I_{\Lambda_L}$ and $I_{\Lambda_L}$ is saturated, so $Q_L\subseteq I_{\Lambda_L}$. Conversely, the proof of \cite[Theorem 3.19]{HHO} shows that the vectors $z$ with $x^u(x^{z^+}-x^{z^-})\in I_L$ for some monomial $x^u$ form a subgroup of $\Z^L$, which contains the generators of $\Lambda_L$; so $I_{\Lambda_L}\subseteq Q_L$ (see also \cite[Corollary 2.5]{ES}). Hence $Q_L=I_{\Lambda_L}=I_{\ker\Phi_L}=\ker\phi_L$, and $\htop Q_L=|L|-1-\kappa(L)$, since $\Phi_L(\Z^L)$ has rank $1+\kappa(L)$.
\end{proof}

Every closed set $W$ is a modular lattice, so Theorem~\ref{thm:toric} applies to it; we write $v_W$ for its Jordan--H\"older valuation.

\begin{proposition}\label{prop:minprimes}
\begin{enumerate}
\item[(a)] For every closed $W$, the ideal $P_W$ is prime, contains $I_L$, and $x_a\in P_W$ if and only if $a\notin W$. In particular $P_W\ne P_{W'}$ for $W\ne W'$.
\item[(b)] $\Min(I_L)$ is the set of minimal elements of $\{P_W : W\subseteq L\text{ closed}\}$.
\item[(c)] $\htop P_W=n-1-\kappa(W)$, and $P_L\in\Min(I_L)$.
\end{enumerate}
\end{proposition}

\begin{proof}
(a) Since $W$ is closed, a generator of $I_L$ either lies in $I_W$ or has both terms divisible by variables outside $W$; so $I_L\subseteq P_W$. By Theorem~\ref{thm:toric}, $S_L/P_W\cong R_W/Q_W$ is a domain in which no variable indexed by $W$ vanishes.

(b) Let $P\supseteq I_L$ be prime and $W=\{a:x_a\notin P\}$. Primality applied to $x_ax_b=x_{a\wedge b}x_{a\vee b}$ makes $W$ closed if it is nonempty, and then $P_W\subseteq P$ by saturation. If $W$ is empty, $P$ is the maximal homogeneous ideal and strictly contains $P_L$. So every minimal prime of $I_L$ is some $P_W$, and the assertion follows.

(c) Theorem~\ref{thm:toric} gives $\dim S_L/P_W=1+\kappa(W)$. No $P_W$ with $W\ne L$ is contained in $P_L$, by (a).
\end{proof}

By Proposition~\ref{prop:minprimes}(b), we have to compare the primes $P_W$. For a finite set $H$ we call \emph{weight} on $\R^H$ an affine function $g(z)=\lambda_0+\sum_{e\in H}\lambda_ez_e$ with rational coefficients.

\begin{proposition}\label{prop:inclusion}
Let $W\subseteq T$ be closed subsets of $L$, and put $h_a=(1,v_T(a))\in\Z\oplus\Z^{\cC(T)}$ for $a\in T$. Then $\dim_\Q\operatorname{span}_\Q\{h_a:a\in W\}\le1+\kappa(W)$. If equality holds, and there is a weight $g$ on $\R^{\cC(T)}$ with $g(v_T(a))\ge0$ for $a\in T$ and $g(v_T(a))=0$ exactly for $a\in W$, then $P_T\subseteq P_W$.
\end{proposition}

\begin{proof}
Let $\Gamma$ be the kernel of $\Z^W\to\Z\oplus\Z^{\cC(T)}$, $e_a\mapsto h_a$. By Theorem~\ref{thm:toric} for $T$, $Q_T\cap R_W$ is the kernel of $x_a\mapsto t\,s^{v_T(a)}$, $a\in W$, so $Q_T\cap R_W=I_\Gamma$. Lemma~\ref{lem:valuation} for $T$ gives $\Lambda_W\subseteq\Gamma$. Since $\rk\Gamma=|W|-\dim\operatorname{span}\{h_a:a\in W\}$ and $\rk\Lambda_W=|W|-1-\kappa(W)$ by Theorem~\ref{thm:toric} for $W$, the inequality follows, and in case of equality $\Gamma=\Lambda_W$, because $\Lambda_W$ is saturated; then $Q_T\cap R_W=Q_W$.

Let $\pi:R_T\to R_W$ send $x_a$ to $0$ for $a\in T\setminus W$ and fix the other variables. The ideal $Q_T$ is generated by binomials $x^\alpha-x^\beta$ with $\sum_a\alpha_ah_a=\sum_a\beta_ah_a$. Since $g(v_T(a))$ is a linear function of $h_a$, we get $\sum_a\alpha_ag(v_T(a))=\sum_a\beta_ag(v_T(a))$; all terms are nonnegative and vanish exactly on $W$, so $x^\alpha$ is supported on $W$ if and only if $x^\beta$ is. Hence $\pi(x^\alpha-x^\beta)$ is $0$ or lies in $Q_T\cap R_W=Q_W$. Since $f-\pi(f)\in(x_a:a\in T\setminus W)$ for $f\in R_T$, we get $Q_T\subseteq P_W$, and therefore $P_T\subseteq P_W$.
\end{proof}

\begin{lemma}\label{lem:maxchain}
	Let $T$ be a finite lattice and $W\subseteq T$ a closed set containing a maximal chain of $T$. Then $W=T$.
\end{lemma}

\begin{proof}
	Let $\hat0=c_0\lessdot c_1\lessdot\cdots\lessdot c_\ell=\hat1$ be the chain. By induction on $\ell$, applied to $W\cap[c_1,\hat1]$, we have $[c_1,\hat1]\subseteq W$. If $z\in T$ is incomparable with $c_1$, then $z\wedge c_1=\hat0$ and $z\vee c_1\in[c_1,\hat1]$ lie in $W$, so $z\in W$; the other elements of $T$ lie in $\{\hat0\}\cup[c_1,\hat1]$.
\end{proof}

\begin{proposition}\label{prop:general}
	Let $W\subseteq L$ be closed with $W\ne L$. Then $\kappa(W)\le\rk W\le r-1$, and $\kappa(W)=r-1$ if and only if $W$ is distributive of rank $r-1$. Consequently $\htop P_W\ge n-r$, and if $L$ is not distributive, then $\htop I_L\ge n-r$.
\end{proposition}

\begin{proof}
	A chain of length $r$ in $L$ is a maximal chain, so $\rk W\le r-1$ by Lemma~\ref{lem:maxchain}, and Lemma~\ref{lem:kappa} gives the rest of the first assertion. Then $\htop P_W\ge n-r$ by Proposition~\ref{prop:minprimes}(c), and also $\htop P_L\ge n-r$ if $L$ is not distributive, by Lemma~\ref{lem:kappa}. Proposition~\ref{prop:minprimes}(b) gives the last assertion.
\end{proof}

\begin{lemma}\label{lem:hibi}
	Let $W\subseteq D\subseteq L$ be closed, with $D$ distributive. Then $P_D\subseteq P_W$. In particular, if $P_W\in\Min(I_L)$ and $W$ is distributive, then $W$ is maximal among the closed distributive subsets of $L$.
\end{lemma}

\begin{proof}
	By Hibi's theorem $I_D$ is prime, and it contains no variable, so $Q_D=I_D$. Let $f=x_ax_b-x_{a\wedge b}x_{a\vee b}$ with $a,b\in D$ incomparable. If $a,b\in W$, then $f\in I_W\subseteq P_W$. Otherwise, since $W$ is closed, also $a\wedge b$ and $a\vee b$ are not both in $W$, so each term of $f$ is divisible by some $x_c$ with $c\in D\setminus W$, and again $f\in P_W$. Hence $P_D\subseteq P_W$. If $W\subsetneq D$, then $P_D\subsetneq P_W$ by Proposition~\ref{prop:minprimes}(a), which gives the second assertion.
\end{proof}

\subsection{Closed subsets of distributive lattices}
We now describe the closed subsets of $\mathcal J(\Pi)$. This will provide the weights needed in Proposition~\ref{prop:inclusion}. For $t\subseteq H$ we write $\one_t\in\{0,1\}^H$ for the indicator vector of $t$.

\begin{lemma}\label{lem:boolean}
Let $T$ be a sublattice of $2^H$. Then $\dim_\Q\operatorname{span}_\Q\{(1,\one_t):t\in T\}=1+\rk T$. If moreover $\emptyset,H\in T$ and $\rk T=|H|$, then $T=\mathcal J(\Pi)$ for a partial order $\Pi$ on $H$.
\end{lemma}

\begin{proof}
Let $t_0\subsetneq\cdots\subsetneq t_k$ be a maximal chain of $T$, so that $k=\rk T$, and put $B_i=t_i\setminus t_{i-1}$. For $t\in T$, the set $(t\cap t_i)\cup t_{i-1}\in T$ lies between $t_{i-1}$ and $t_i$, so $t$ contains all or none of $B_i$. Hence each $\one_t-\one_{t_0}$ is a sum of some of the vectors $\one_{B_i}$, which are nonzero with disjoint supports, and the first assertion follows. For the second, put $H_j=\bigcap\{t\in T:j\in t\}\in T$, and let $i\le_\Pi j$ if $i\in H_j$. This relation is reflexive and transitive. A maximal chain of $T$ goes from $\emptyset$ to $H$ adding one element at a time; for $i\ne j$, its first member containing $i$ or $j$ contains only one of them, say $j$, so $i\notin H_j$. Thus $\le_\Pi$ is a partial order. If $j\in t\in T$, then $H_j\subseteq t$; so the members of $T$ are order ideals, and every nonempty order ideal $A$ equals $\bigcup_{j\in A}H_j\in T$.
\end{proof}

\begin{lemma}\label{lem:closed-distributive}
Let $\Pi$ be a poset, and let $\widehat\Pi$ be obtained by adding a least element $\bot$ and a greatest element $\top$. For $A\in\mathcal J(\Pi)$ put $z_h(A)=\one_A(h)$ for $h\in\Pi$, $z_\bot(A)=1$ and $z_\top(A)=0$. If $U\subseteq\mathcal J(\Pi)$ is closed, there is a set $F$ of pairs $h\lessdot k$ of $\widehat\Pi$ such that $U$ consists of the order ideals $A$ with $z_h(A)=z_k(A)$ for all $(h,k)\in F$. Thus $U$ is the zero set of the weight $\sum_{(h,k)\in F}(z_h-z_k)$, which is nonnegative on $\mathcal J(\Pi)$. If moreover $U$ is maximal among the closed subsets of $\mathcal J(\Pi)$ not containing some $A_*$, then $\rk U=|\Pi|-1$.
\end{lemma}

\begin{proof}
Let $A_0\subsetneq\cdots\subsetneq A_s$ be a maximal chain of $U$ and $B_i=A_i\setminus A_{i-1}$; as above, every member of $U$ contains all or none of each $B_i$. Each $B_i$ is connected in the comparability graph of $\Pi$: otherwise $B_i=B'\sqcup B''$ with no element of $B'$ comparable with one of $B''$, and $A_{i-1}\cup B'$, $A_{i-1}\cup B''$ are incomparable order ideals with intersection $A_{i-1}$ and union $A_i$, so they lie in $U$, against the maximality of the chain. The order ideals between $A_0$ and $A_s$ containing all or none of each $B_i$ form a sublattice $T\supseteq U$ in which our chain is maximal, and $U$ is closed in $T$, so $U=T$ by Lemma~\ref{lem:maxchain}. Hence $U$ is defined by the equations $z_h=z_\bot$ for $h\in A_0$, $z_h=z_\top$ for $h\notin A_s$, and $z_h=z_k$ for $h,k$ in the same $B_i$, and by connectedness the last ones may be taken along comparable pairs. Since $z_h(A)\ge z_k(A)$ for $h\le k$, an equation along a comparable pair is equivalent to the equations along the covers of a saturated chain between them. This gives $F$.

Let now $U$ be as in the last assertion. Some equation $z_h=z_k$ of $F$ fails at $A_*$, and its zero set $Z$ contains $U$ and not $A_*$. Since $z_h-z_k\ge0$ and $z_h(A)+z_h(A')=z_h(A\cap A')+z_h(A\cup A')$, the set $Z$ is closed, so $Z=U$. If $h=\bot$, then $k$ is minimal in $\Pi$ and $Z\cong\mathcal J(\Pi\setminus\{k\})$; if $k=\top$, then $Z=\mathcal J(\Pi\setminus\{h\})$. If $h,k\in\Pi$, then $I=\{x\in\Pi:x<k\}\setminus\{h\}$ and ${\downarrow}k=I\cup\{h,k\}$ are order ideals, and refining $\emptyset\subseteq I$ and ${\downarrow}k\subseteq\Pi$ one element at a time gives a chain of $Z$ of length $|\Pi|-1$, which is maximal because no member of $Z$ lies strictly between $I$ and ${\downarrow}k$. In all cases $\rk Z=|\Pi|-1$.
\end{proof}

\section{Lattices with one diamond}\label{sec:structure}

From now on $L$ is a finite modular lattice with exactly one diamond $[p,q]$, $A_\dm$ is the set of atoms of $[p,q]$, so that $|A_\dm|\ge3$, and $n=|L|$, $r=\rk L$. For a map $f$ on $L$, we call $f^{-1}(y)$ the \emph{fibre} of $f$ over $y$. In this section we embed $L$ into $[p,q]\times\Delta$, for a distributive lattice $\Delta$ (Lemma~\ref{lem:rho}); this embedding is the main tool of the paper.

\begin{lemma}\label{lem:onediamond}
	All prime intervals of $[p,q]$ lie in one projectivity class $C_\dm$. Every maximal chain of $L$ contains exactly two prime intervals of class $C_\dm$ and exactly one of every other class. In particular $\kappa(L)=r-1$.
\end{lemma}

\begin{proof}
	For pairwise distinct $a,b,c\in A_\dm$ we have $[p,a]\nearrow[b,q]\searrow[p,c]\nearrow[a,q]$, so the prime intervals of $[p,q]$ form one class $C_\dm$, and a maximal chain through $p\lessdot a\lessdot q$ shows $m_{C_\dm}(\hat0,\hat1)\ge2$. By Lemma~\ref{lem:valuation}, $\kappa(L)\le r-1$. For the reverse inequality we use a result of Herrmann \cite{He}. A \emph{line} of $L$ is a maximal set of at least three join-irreducible elements with pairwise equal joins, and this common join is a \emph{line-top}. By \cite[Corollary 1.2]{He}, the number $i(L)$ of line-tops satisfies $i(L)\ge r-s(L)$, where $s(L)$ is the number of maximal proper congruences of $L$. If $t$ is a line-top and $t_*$ is the meet of the elements covered by $t$, then $[t_*,t]$ has length two and at least three atoms \cite[Section 2, assertion (2)]{He}, so $i(L)\le1$. The join-irreducible congruences of a finite lattice are the congruences generated by prime intervals, and two prime intervals generate the same congruence if and only if each is congruence-projective to the other \cite[Lemma 238]{G}; in a finite modular lattice this means that they are projective, and the congruence lattice is Boolean \cite[Section V.1.5]{G}. Hence $s(L)$, the number of coatoms of the congruence lattice, equals the number of its atoms, which is $\kappa(L)$, and $\kappa(L)\ge r-1$. Finally, the $r-1$ multiplicities $m_C(\hat0,\hat1)$ are positive, add up to $r$, and $m_{C_\dm}(\hat0,\hat1)\ge2$; so $m_{C_\dm}(\hat0,\hat1)=2$ and the others equal $1$.
\end{proof}

Let $\cE=\cC(L)\setminus\{C_\dm\}$, so that $|\cE|=r-2$, and write $v=(v_\dm,(v_E)_{E\in\cE})$ with $v_\dm=v_{C_\dm}$. We call $v_\dm(x)$ the \emph{level} of $x$. By Lemma~\ref{lem:onediamond}, $v_\dm(x)\in\{0,1,2\}$ and $v_E(x)\in\{0,1\}$, and we put
\[
\beta(x)=\{E\in\cE:v_E(x)=1\},\qquad t_0=\beta(p),
\]
and $L_i=\{x\in L:v_\dm(x)=i\}$ for $i=0,1,2$. The level and $\beta$ are monotone, since $v(y)-v(x)=m(x,y)\ge0$ for $x\le y$.

\begin{lemma}\label{lem:M3}
	\begin{enumerate}
		\item[(a)] Every sublattice of $L$ isomorphic to $M_3$ has least element $p$, greatest element $q$, and its other three elements in $A_\dm$.
		\item[(b)] For a closed set $W$ the following are equivalent: $W$ is not distributive; $p,q\in W$; $|W\cap A_\dm|\ge2$; $[p,q]\subseteq W$.
		\item[(c)] If $W$ is closed and $[p,q]\subseteq W$, then $[p,q]$ is the only diamond of the lattice $W$, and $\kappa(W)=\rk W-1$.
	\end{enumerate}
\end{lemma}

\begin{proof}
	(a) Let $\{o,x,y,z,i\}\cong M_3$, with least element $o$ and greatest element $i$. By Lemma~\ref{lem:valuation}, $m(o,a)=m(b,i)$ for distinct $a,b\in\{x,y,z\}$, so the six vectors $m(o,a)$ and $m(a,i)$ equal a common vector $u\ne0$, and $m(o,i)=2u$. The entries of $m(o,i)$ are bounded by those of $m(\hat0,\hat1)$, which are $2$ at $C_\dm$ and $1$ elsewhere; hence $u=e_{C_\dm}$, and $[o,i]$ is an interval of length two with atoms $x,y,z$. So $[o,i]=[p,q]$.

	(b) If $W$ is not distributive, it contains a copy of $M_3$ by Birkhoff's theorem, and hence $p$ and $q$ by (a). If $p,q\in W$, closedness gives $[p,q]\subseteq W$, since distinct atoms of the diamond have meet $p$ and join $q$. If $[p,q]\subseteq W$, then $|W\cap A_\dm|\ge2$. Finally, two atoms in $W$ give $p,q\in W$, hence $[p,q]\subseteq W$, and $W$ contains a copy of $M_3$.

	(c) The atoms of a diamond of $W$ have pairwise meet and join equal to its endpoints, in $W$ and hence in $L$; by (a) the diamond is $[p,q]$. Lemma~\ref{lem:onediamond} for $W$ gives $\kappa(W)=\rk W-1$.
\end{proof}

\begin{lemma}\label{lem:beta}
	\begin{enumerate}
		\item[(a)] $\beta(x\wedge y)=\beta(x)\cap\beta(y)$ and $\beta(x\vee y)=\beta(x)\cup\beta(y)$ for all $x,y\in L$.
		\item[(b)] If $x\le y$, the length of $[x,y]$ is $v_\dm(y)-v_\dm(x)+|\beta(y)\setminus\beta(x)|$.
		\item[(c)] $v_\dm(p)=0$, $v_\dm(a)=1$ for $a\in A_\dm$, $v_\dm(q)=2$, and $\beta(z)=t_0$ for $z\in[p,q]$.
	\end{enumerate}
\end{lemma}

\begin{proof}
	(a) Fix $E\in\cE$. If $v_E(x)=v_E(y)=1$, then $v_E(x\vee y)=1$, and $v_E(x\wedge y)=v_E(x)+v_E(y)-v_E(x\vee y)=1$ by Lemma~\ref{lem:valuation}; otherwise $v_E(x\wedge y)=0$ by monotonicity. So $v_E(x\wedge y)=\min\{v_E(x),v_E(y)\}$, and dually $v_E(x\vee y)=\max\{v_E(x),v_E(y)\}$. (b) The length is the sum of the entries of $m(x,y)=v(y)-v(x)$. (c) The only prime intervals of class $C_\dm$ on a maximal chain through $p\lessdot a\lessdot q$ are $[p,a]$ and $[a,q]$, so $v_\dm(p)=0$, $v_\dm(a)=1$, $v_\dm(q)=2$, and $m(p,q)=2e_{C_\dm}$ gives $\beta(q)=\beta(p)$; monotonicity does the rest.
\end{proof}

\begin{lemma}\label{lem:levels}
	\begin{enumerate}
		\item[(a)] $L_0={\downarrow}p$ and $L_2={\uparrow}q$.
		\item[(b)] $\beta$ is injective on ${\downarrow}p$ and on ${\uparrow}q$.
		\item[(c)] For $t\in\beta(L)$, the fibre $\beta^{-1}(t)$ is an interval. Moreover $\beta^{-1}(t_0)=[p,q]$, and if $t\ne t_0$, then $\beta^{-1}(t)$ is a single element or a pair $m\lessdot M$ with $v_\dm(M)=v_\dm(m)+1$.
	\end{enumerate}
\end{lemma}

\begin{proof}
	(a) Clearly ${\downarrow}p\subseteq L_0$. Let $x\in L_0$, and put $s=x\vee p$, $S=x\vee q$ and $s_a=x\vee a$ for $a\in A_\dm$. Since $x\wedge p$, $x\wedge a$ and $x\wedge q$ have level $0$, Lemmas~\ref{lem:valuation} and~\ref{lem:beta} show that $s$, $s_a$, $S$ have levels $0,1,2$ and the same $\beta$-value $\beta(x)\cup t_0$. By Lemma~\ref{lem:beta}(b), $[s,S]$ has length two and each $s_a$ is an atom of it; the $s_a$ are distinct, since $s_a=s_b$ with $a\ne b$ would lie above $q$ and have level $2$. So $[s,S]=[p,q]$, and $x\le s=p$. The case of $L_2$ is dual.

	(b) If $x,x'\le p$ and $\beta(x)=\beta(x')$, then $\beta(x\wedge x')=\beta(x)$, and $[x\wedge x',x]$ has length zero by Lemma~\ref{lem:beta}(b); so $x=x\wedge x'=x'$. The case of ${\uparrow}q$ is dual.

	(c) By Lemma~\ref{lem:beta}(a) the fibre is a convex sublattice, hence an interval $[m,M]$, and its length is $v_\dm(M)-v_\dm(m)\le2$. If the length is two, then $m\le p\le q\le M$ by (a), so $t=\beta(m)\subseteq t_0\subseteq\beta(M)=t$, and (b) gives $m=p$, $M=q$. Since $\beta$ is constant on $[p,q]$, the fibre of $t_0$ contains $[p,q]$, so it has length two and equals $[p,q]$; the other fibres have length at most one.
\end{proof}

Let $\Delta=\beta(L)\subseteq2^\cE$, and write $[s,s']_\Delta=\{t\in\Delta:s\subseteq t\subseteq s'\}$ for $s\subseteq s'$ in $\Delta$.

\begin{lemma}\label{lem:Delta}
\begin{enumerate}
\item[(a)] $\Delta$ is a sublattice of $2^\cE$ containing $\emptyset$ and $\cE$, every cover $s\lessdot s'$ of $\Delta$ has $|s'\setminus s|=1$, and $\Delta=\mathcal J(\Pi_\Delta)$ for a partial order $\Pi_\Delta$ on $\cE$.
\item[(b)] If $U\subseteq\Delta$ is closed, there is a weight $f$ on $\R^\cE$, a sum of terms of the forms $z_i-z_j$, $1-z_j$ and $z_i$, with $f(\one_s)\ge0$ for all $s\in\Delta$ and $f(\one_s)=0$ if and only if $s\in U$. If $U$ is maximal among the closed subsets of $\Delta$ not containing some $s_*\in\Delta$, then $\rk U=r-3$.
\item[(c)] If $f$ is a weight on $\R^\cE$ with $f(\one_s)\ge0$ for all $s\in\Delta$, then $\{s\in\Delta:f(\one_s)=0\}$ is empty or closed.
\item[(d)] If $U$ is closed in $\Delta$, then $\beta^{-1}(U)$ is closed in $L$. If $W$ is closed in $L$ and $W=\beta^{-1}(\beta(W))$, then $\beta(W)$ is closed in $\Delta$.
\end{enumerate}
\end{lemma}

\begin{proof}
(a) By Lemma~\ref{lem:beta}(a), $\Delta$ is a sublattice, and it contains $\beta(\hat0)=\emptyset$ and $\beta(\hat1)=\cE$. If $s\lessdot s'$ in $\Delta$, choose $x\le y$ (replacing $x$ by $x\wedge y$) with $\beta(x)=s$ and $\beta(y)=s'$. Along a maximal chain of $[x,y]$, the $\beta$-values increase by at most one element at each step, and each of them is $s$ or $s'$; so $|s'\setminus s|=1$. Hence $\rk\Delta=|\cE|$, and Lemma~\ref{lem:boolean} gives $\Pi_\Delta$. (b) This is Lemma~\ref{lem:closed-distributive} for $\mathcal J(\Pi_\Delta)$, since $|\Pi_\Delta|=r-2$. (c) The values of $f$ are nonnegative and satisfy $f(\one_s)+f(\one_{s'})=f(\one_{s\cap s'})+f(\one_{s\cup s'})$, so $f$ vanishes at $s,s'$ if and only if it vanishes at $s\cap s',s\cup s'$.

(d) Let $x,y\in L$ be incomparable. If $\beta(x)\subseteq\beta(y)$, then $\beta(x\wedge y)=\beta(x)$ and $\beta(x\vee y)=\beta(y)$; otherwise $\beta(x),\beta(y)$ are incomparable, with meet $\beta(x\wedge y)$ and join $\beta(x\vee y)$. In both cases closedness of $U$ gives $x,y\in\beta^{-1}(U)$ if and only if $x\wedge y,x\vee y\in\beta^{-1}(U)$. For the second assertion, $\beta(W)$ is a sublattice; if $s,s'\in\Delta$ are incomparable with $s\cap s',s\cup s'\in\beta(W)$, choose $x,y$ with $\beta(x)=s$ and $\beta(y)=s'$. Then $x,y$ are incomparable and $x\wedge y,x\vee y\in\beta^{-1}(\beta(W))=W$, so $x,y\in W$.
\end{proof}

\begin{lemma}\label{lem:rho}
	The map $\rho:L\to[p,q]$, $\rho(x)=(x\vee p)\wedge q$, is a lattice homomorphism with $\rho^{-1}(p)=L_0$, $\rho^{-1}(q)=L_2$ and $\rho^{-1}(A_\dm)=L_1$, and $(\rho,\beta):L\to[p,q]\times\Delta$ is an injective lattice homomorphism. For $c\in A_\dm$, the fibre $\rho^{-1}(c)$ is an interval $[m_c,M_c]$ containing $c$, and $\lt_c=\beta(m_c)$, $\ut_c=\beta(M_c)$ satisfy $\lt_c\subseteq t_0\subseteq\ut_c$ and $\beta(\rho^{-1}(c))=[\lt_c,\ut_c]_\Delta$. If $c\ne c'$, then $\lt_c\cup\lt_{c'}=t_0$ and $\ut_c\cap\ut_{c'}=t_0$.
\end{lemma}

\begin{proof}
	Write $d=v_\dm$. By Lemma~\ref{lem:levels}(a), $\rho=p$ on $L_0$ and $\rho=q$ on $L_2$. For $x\in L_1$, Lemma~\ref{lem:valuation} gives $d(x\vee p)=1+0-0=1$ and $d(\rho(x))=1+2-2=1$, so $\rho(x)\in A_\dm$. This gives the fibres of $p$, $q$ and $A_\dm$.

	Let $x,y$ be incomparable; we check that $\rho$ preserves their meet and join. If both lie in $L_0$, or both in $L_2$, this is clear, and one in $L_0$ and one in $L_2$ is impossible. If $x\in L_0$ and $y\in L_1$, then $x\wedge y\in L_0$ and $d(x\vee y)=1$, so $\rho(x\vee y)$ is an atom above $\rho(y)$, hence equal to it; the case $x\in L_1$, $y\in L_2$ is dual. Let $x,y\in L_1$. If $\rho(x)\ne\rho(y)$, then $\rho(x\vee y)\ge\rho(x)\vee\rho(y)=q$, so $d(x\vee y)=2$, $d(x\wedge y)=0$ and $\rho(x\wedge y)=p$. If $\rho(x)=\rho(y)=a$, then $x\vee p$, $y\vee p$ and $q$ lie in $[a,\hat1]$, which is distributive because it contains no copy of $M_3$ by Lemma~\ref{lem:M3}(a). Hence $\rho(x\vee y)=((x\vee p)\vee(y\vee p))\wedge q=\rho(x)\vee\rho(y)=a$, so $d(x\wedge y)=1$ and $\rho(x\wedge y)=a$.

	If $\rho(x)=\rho(y)$ and $\beta(x)=\beta(y)=t$, then $x,y\in[p,q]$ if $t=t_0$, where $\rho$ is the identity; otherwise $x,y$ have the same level, which is determined by $\rho$, and $x=y$ by Lemma~\ref{lem:levels}(c). This proves injectivity.

	The fibre $\rho^{-1}(c)$ is a convex sublattice containing $c$, hence an interval $[m_c,M_c]$, and $\beta(\rho^{-1}(c))\subseteq[\lt_c,\ut_c]_\Delta$. Conversely, if $t\in[\lt_c,\ut_c]_\Delta$ and $\beta(z)=t$, then $(z\vee m_c)\wedge M_c$ lies in the fibre and has $\beta$-value $t$. If $c\ne c'$, two elements of level one in the fibres of $c$ and $c'$ with the same $\beta$-value $t$ force $t=t_0$ by Lemma~\ref{lem:levels}(c); so $[\lt_c,\ut_c]_\Delta\cap[\lt_{c'},\ut_{c'}]_\Delta=\{t_0\}$. This intersection contains $\lt_c\cup\lt_{c'}$ and $\ut_c\cap\ut_{c'}$, which are therefore equal to $t_0$.
\end{proof}

For $c\in A_\dm$ we put $V_c=\ut_c\setminus t_0$ and $R_c=t_0\setminus\lt_c$. By Lemma~\ref{lem:rho}, the sets $R_c\cup V_c$, $c\in A_\dm$, are pairwise disjoint, and $\beta(x)\setminus t_0\subseteq V_c$ and $t_0\setminus\beta(x)\subseteq R_c$ whenever $\rho(x)=c$.

\section{Inclusions among the candidate primes}\label{sec:inclusions}

By Lemma~\ref{lem:M3}(b), a closed set $W\ne L$ either contains the diamond, or it is distributive and contains at most one atom of the diamond and not both $p$ and $q$. In this section we decide when $P_L\subseteq P_W$, using Proposition~\ref{prop:inclusion} with $T=L$. We write the coordinates of $\R^{\cC(L)}$ as $z_\dm$ and $z_E$, $E\in\cE$. If a weight $g$ satisfies $g(v(x))\ge0$ for all $x\in L$, then $\{x:g(v(x))=0\}$ is empty or closed, since $g(v(x))+g(v(y))=g(v(x\wedge y))+g(v(x\vee y))$ by Lemma~\ref{lem:valuation}.

\subsection{Closed sets containing the diamond}

\begin{lemma}\label{lem:saturation}
If $W$ is closed and $[p,q]\subseteq W$, then $W=\beta^{-1}(\beta(W))$.
\end{lemma}

\begin{proof}
Let $w\in W$ and $\beta(z)=\beta(w)=t$. If $t=t_0$, then $z\in[p,q]$. Otherwise, by Lemma~\ref{lem:levels}(c), we must show that $m\in W$ if and only if $M\in W$ for a pair $m\lessdot M$ with $\beta(m)=\beta(M)=t$; by duality, let $m\in L_0$ and $M\in L_1$, so that $t\subseteq t_0$. If $M\in W$, then $M\wedge p\in W$ has $\beta$-value $t$, so $M\wedge p=m$ by Lemma~\ref{lem:levels}(b). If $m\in W$, take $a\in A_\dm\setminus\{\rho(M)\}$; then $(\rho,\beta)(M\wedge a)=(p,t)$ and $(\rho,\beta)(M\vee a)=(q,t_0)$, so $M\wedge a=m$ and $M\vee a=q$; as $M$ and $a$ are incomparable, $M\in W$ by closedness.
\end{proof}

\begin{proposition}\label{prop:diamondinclusion}
If $W$ is closed and $[p,q]\subseteq W$, then $P_L\subseteq P_W$.
\end{proposition}

\begin{proof}
By Lemma~\ref{lem:saturation} and Lemma~\ref{lem:Delta}(d),(b), there is a weight $f$ on $\R^\cE$, nonnegative on $\Delta$, whose zero set in $\Delta$ is $\beta(W)$. Then $g(z)=f((z_E)_{E\in\cE})$ satisfies $g(v(x))=f(\one_{\beta(x)})$, so it is nonnegative on $L$ and vanishes exactly on $\beta^{-1}(\beta(W))=W$. For the rank condition, take a maximal chain of the lattice $W$ through $p\lessdot a\lessdot q$; its length is $\rk W$. The two steps $p\lessdot a\lessdot q$ have increment $e_{C_\dm}$ under $v$; the other steps lie in ${\downarrow}p$ or ${\uparrow}q$, so their increments have $C_\dm$-entry zero and nonempty supports in $\cE$, which are pairwise disjoint because each $v_E$ increases at most once. Hence the increments span a space of dimension $\rk W-1=\kappa(W)$, by Lemma~\ref{lem:M3}(c), and together with $(1,v(\min W))$ they show that the vectors $(1,v(w))$, $w\in W$, span a space of dimension at least $1+\kappa(W)$. Proposition~\ref{prop:inclusion} gives $P_L\subseteq P_W$.
\end{proof}

\begin{lemma}\label{lem:oneatom}
Let $W\subseteq W'$ be closed, with $|W\cap A_\dm|=1$ and $|W'\cap A_\dm|\ge2$. Then $P_{W'}\not\subseteq P_W$.
\end{lemma}

\begin{proof}
Let $W\cap A_\dm=\{a\}$, and choose distinct $b,c\in A_\dm\setminus\{a\}$. By Lemma~\ref{lem:M3}(b), $[p,q]\subseteq W'$, and in the lattice $W'$ both $[p,a]$ and $[p,b]$ transpose up to $[c,q]$. So $v_{W'}(a)=v_{W'}(b)$, and $x_a-x_b\in Q_{W'}\subseteq P_{W'}$ by Theorem~\ref{thm:toric}. But $x_b\in P_W$ and $x_a\notin P_W$ by Proposition~\ref{prop:minprimes}(a), so $x_a-x_b\notin P_W$.
\end{proof}

\subsection{Closed sets avoiding the diamond}
Let $W$ be closed with $W\cap[p,q]=\emptyset$, hence distributive. To apply Proposition~\ref{prop:inclusion} we separate $\beta(W)$ from $t_0$ by a weight on $\Delta$ (Lemma~\ref{lem:separation}), starting with the case in which $W$ does not meet $L_1$.

\begin{lemma}\label{lem:outer}
Let $W$ be closed with $W\cap[p,q]=\emptyset$ and $W\subseteq L_0\cup L_2$. Then there is a weight $f$ on $\R^\cE$ with $f(\one_s)\ge0$ for $s\in\Delta$, $f(\one_{\beta(w)})=0$ for $w\in W$, and $f(\one_{t_0})=1$.
\end{lemma}

\begin{proof}
If $W\subseteq L_0$, then $u=\max W<p$, so $\beta(u)\subsetneq t_0$ by Lemma~\ref{lem:levels}(b), and $f=z_i$ with $i\in t_0\setminus\beta(u)$ works. If $W\subseteq L_2$, take dually $f=1-z_j$ with $j\in\beta(\min W)\setminus t_0$. Otherwise let $u=\max(W\cap L_0)$, $w=\min(W\cap L_2)$, $s=\beta(u)$, $t=\beta(w)$, $R=t_0\setminus s$ and $S=t\setminus t_0$; as before, $R,S\ne\emptyset$. We claim that there are $i\in R$ and $j\in S$ such that every element of $\Delta$ containing $j$ contains $i$; then $f=z_i-z_j$ works, since $\beta(w')\subseteq s$ for $w'\in W\cap L_0$ and $\beta(w')\supseteq t$ for $w'\in W\cap L_2$. If not, then for each $j\in S$ the set $H_j=\bigcap\{h\in\Delta:j\in h\}\in\Delta$ is contained in $t$ and disjoint from $R$, so $h=s\cup\bigcup_{j\in S}H_j\in\Delta$ satisfies $h\cap t_0=s$ and $h\setminus t_0=S$. If $\beta(x)=h$, then $x\in L_1$, because $h$ neither contains nor is contained in $t_0$, and for $b\in A_\dm\setminus\{\rho(x)\}$ we get $(\rho,\beta)(x\wedge b)=(p,s)$ and $(\rho,\beta)(x\vee b)=(q,t)$, that is, $x\wedge b=u$ and $x\vee b=w$. Closedness gives $b\in W$, a contradiction.
\end{proof}

\begin{lemma}\label{lem:separation}
Let $W$ be closed with $W\cap[p,q]=\emptyset$. Then there is a weight $f$ as in Lemma~\ref{lem:outer}. In particular, the smallest closed subset of $\Delta$ containing $\beta(W)$ does not contain $t_0$.
\end{lemma}

\begin{proof}
The lattice $\rho(W)$ is distributive, so it contains at most two atoms of the diamond, since three atoms generate a copy of $M_3$; say $\rho(W)\subseteq\{p,a,b,q\}$ with $a\ne b$. Let $B=\rho^{-1}(\{p,a,b,q\})$, and identify $\{p,a,b,q\}$ with $2^{\{\sigma,\tau\}}$ via $a\mapsto\{\sigma\}$, $b\mapsto\{\tau\}$. Then $(\rho,\beta)$ gives an injective lattice homomorphism $\gamma:B\to2^H$, $H=\cE\sqcup\{\sigma,\tau\}$, with $\gamma(x)\cap\cE=\beta(x)$. Its image contains $\emptyset$ and $H$, and $\rk B=r=|H|$, since $B$ contains a maximal chain of $L$ through $p\lessdot a\lessdot q$. By Lemma~\ref{lem:boolean}, $\gamma(B)=\mathcal J(\Pi)$ for a partial order $\Pi$ on $H$, in which $\sigma$ and $\tau$ are incomparable, as shown by $\gamma(a)$ and $\gamma(b)$.

The set $\gamma(W)$ is closed in $\mathcal J(\Pi)$, and Lemma~\ref{lem:closed-distributive} describes it by equations $z_h=1$, $z_k=0$ and $z_h=z_k$ with $h\lessdot k$ in $\Pi$. We claim that these equations imply, on $\gamma(W)$, one of the following: (i) $z_i=0$ for some $i\in t_0$; (ii) $z_j=1$ for some $j\in\cE\setminus t_0$; (iii) $z_i=z_j$ for some $i\in t_0$ and $j\in\cE\setminus t_0$ with $i<_\Pi j$. Note that $\gamma(p),\gamma(a),\gamma(b),\gamma(q)$ are the four order ideals $t_0\cup\epsilon$, $\epsilon\subseteq\{\sigma,\tau\}$, and none of them lies in $\gamma(W)$. If an equation involving only coordinates in $\cE$ fails at $\one_{t_0}$, we are in case (i), (ii) or (iii), the orientation in (iii) coming from the fact that $t_0$ is an order ideal of $\Pi$. Otherwise we substitute $\one_{t_0}$ for the coordinates in $\cE$. No equation involves both $\sigma$ and $\tau$, so each remaining equation prescribes $0$ or $1$ to $z_\sigma$ or to $z_\tau$, and since the four ideals above are excluded, some $z_k$ with $k\in\{\sigma,\tau\}$ receives both values. As $\gamma(W)\ne\emptyset$, the two prescriptions are not $z_k=0$ and $z_k=1$. If they are $z_k=0$ and $z_k=z_i$ with $i\in t_0$, we are in case (i), and $z_k=1$, $z_k=z_j$ with $j\notin t_0$ give case (ii). Otherwise they are $z_k=z_i$ and $z_k=z_j$ with $i\in t_0$ and $j\notin t_0$, and evaluating at $\gamma(p)$ and $\gamma(q)$ shows $i<_\Pi k<_\Pi j$; this is case (iii).

In cases (i) and (ii) take $f=z_i$ and $f=1-z_j$. In case (iii), $z_i-z_j$ is nonnegative on $\gamma(B)$ and vanishes on $\gamma(W)$. If it is nonnegative on $\Delta$, take $f=z_i-z_j$. Otherwise some $x\in L$ has $j\in\beta(x)$ and $i\notin\beta(x)$; then $x\in L_1$, since $i\in t_0$ and $j\notin t_0$, and $c=\rho(x)\notin\{a,b\}$, since $z_i-z_j\ge0$ on $\gamma(B)$. Since $i\notin\lt_c$ and $j\in\ut_c$, Lemma~\ref{lem:rho} gives $i\in\lt_a\cap\lt_b$ and $j\notin\ut_a\cup\ut_b$, so $z_i-z_j=1$ on the fibres of $a$ and $b$, which therefore do not meet $W$. Thus $W\subseteq L_0\cup L_2$, and Lemma~\ref{lem:outer} applies. Finally, the zero set of $f$ in $\Delta$ is closed by Lemma~\ref{lem:Delta}(c), contains $\beta(W)$, and does not contain $t_0$.
\end{proof}

\begin{proposition}\label{prop:avoid}
Let $W$ be maximal among the closed distributive subsets of $L$, with $W\cap[p,q]=\emptyset$. Then $W=\beta^{-1}(\beta(W))$ and $\beta(W)$ is closed in $\Delta$. If $\rk W\le r-2$, then $P_L\subseteq P_W$.
\end{proposition}

\begin{proof}
Let $T$ be the smallest closed subset of $\Delta$ containing $\beta(W)$. By Lemmas~\ref{lem:separation}, \ref{lem:Delta}(d) and~\ref{lem:M3}(b), $\beta^{-1}(T)$ is closed, disjoint from $[p,q]=\beta^{-1}(t_0)$, and distributive, so $W=\beta^{-1}(T)$ by maximality, and $\beta(W)=T$. As in Proposition~\ref{prop:diamondinclusion}, Lemma~\ref{lem:Delta}(b) gives a weight, nonnegative on $L$, whose zero set is $W$.

For the rank condition, recall that $\kappa(W)=\rk W$ and that every class of $W$ meets every maximal chain of $W$ once (Lemma~\ref{lem:kappa}). By Lemma~\ref{lem:valuation}, transposed covers of $W$ have the same increment under $v$, so the increment $\mu_C$ depends only on the class $C\in\cC(W)$. Write $\mu_C=d_Ce_{C_\dm}+\one_{B_C}$ with $B_C\subseteq\cE$. On a maximal chain of $W$ each $\mu_C$ occurs once, so the nonempty sets $B_C$ are pairwise disjoint and $\sum_Cd_C\le2$; as each $\mu_C\ne0$, the vectors $\mu_C$ are linearly independent unless two of them equal $e_{C_\dm}$. In that case $\beta$ is constant on exactly two steps of a maximal chain $w_0\lessdot\cdots\lessdot w_s$ of $W$, and the distinct sets among the $\beta(w_i)$ form a maximal chain of $T$: if $t\in T$ lay strictly between $\beta(w_{i-1})$ and $\beta(w_i)$, with $t=\beta(z)$, $z\in W$, then $(z\vee w_{i-1})\wedge w_i$ would lie strictly between $w_{i-1}$ and $w_i$. So $\rk T=\rk W-2$. Moreover $T$ is maximal among the closed subsets of $\Delta$ not containing $t_0$, by the maximality of $W$, so $\rk T=r-3$ by Lemma~\ref{lem:Delta}(b), and $\rk W=r-1$. Hence, if $\rk W\le r-2$, the $\rk W$ vectors $\mu_C$ are independent, and together with $(1,v(\min W))$ they show that the vectors $(1,v(w))$, $w\in W$, span a space of dimension $1+\kappa(W)$. Proposition~\ref{prop:inclusion} gives $P_L\subseteq P_W$.
\end{proof}

\subsection{Closed sets containing an endpoint of the diamond}
Let $W$ be maximal among the closed distributive sets, with $W\cap A_\dm=\emptyset$ and $W\cap\{p,q\}\ne\emptyset$. By duality we may assume $p\in W$, and then $q\notin W$.

\begin{lemma}\label{lem:alphalevels}
Let $W$ be maximal among the closed distributive subsets of $L$, with $p\in W$ and $W\cap A_\dm=\emptyset$. Then $W$ meets both $L_1$ and $L_2$.
\end{lemma}

\begin{proof}
If $W\cap L_2=\emptyset$, the interval ${\downarrow}\max W$ is closed, does not contain $q$, hence is distributive, and it contains $W$. If $\max W\in L_0$, then $\max W=p$ and ${\downarrow}a$ strictly contains $W$ for $a\in A_\dm$; if $\max W\in L_1$, then ${\downarrow}\max W$ contains the atom $\rho(\max W)$. Both contradict maximality.

Suppose that $W\cap L_1=\emptyset$, and let $u=\min(W\cap L_2)$, $t=\beta(u)$ and $S=t\setminus t_0$, which is nonempty since $u\ne q$. If $S\subseteq V_a\cup V_b$ for distinct atoms $a,b$, lift $t\cap\ut_a$ and $t\cap\ut_b$ to $y_a,y_b$ in the fibres of $a$ and $b$ (Lemma~\ref{lem:rho}); then $y_a\wedge y_b=p$ and $y_a\vee y_b=u$, and closedness gives $y_a\in W\cap L_1$. If some $j\in S$ lies in no $V_c$, every $x$ with $j\in\beta(x)$ lies in $L_2$. For $a\in A_\dm$, the set $T=\rho^{-1}(\{p,a\})\cup\{x:j\in\beta(x)\}$ is the union of a sublattice closed downwards and a sublattice closed upwards, hence a sublattice. It is closed: if $x,y$ are incomparable with $x\wedge y\in\rho^{-1}(\{p,a\})$ and $j\in\beta(x\vee y)$, say $j\in\beta(x)$, then $\rho(x)=q$ and $\rho(y)=\rho(x\wedge y)\in\{p,a\}$, and the other cases are clear. Since $\rho(T)$ is a chain, $(\rho,\beta)$ embeds $T$ into a distributive lattice. Finally $T$ contains $a$ and $W$, because the elements of $W\cap L_2$ lie above $u$, a contradiction. In the remaining case $S$ meets at least three of the sets $V_c$; take $i\in S\cap V_a$ and $j\in S\cap V_b$ with $a\ne b$. The weight $z_\dm-z_i-z_j$ is nonnegative on $L$, because at most one of $i,j$ lies in $\beta(x)$ when $x\in L_1$. Its zero set contains $W$ and $M_a$, but not $q$, so it is a closed distributive set strictly containing $W$, a contradiction.
\end{proof}

\begin{proposition}\label{prop:endpoint}
Let $W$ be maximal among the closed distributive subsets of $L$, with $W\cap A_\dm=\emptyset$ and $W\cap\{p,q\}\ne\emptyset$. Then $P_L\subseteq P_W$.
\end{proposition}

\begin{proof}
Let $p\in W$. By Lemma~\ref{lem:alphalevels}, $W$ meets $L_1$ and $L_2$, and since $\rho(W)$ is distributive, $W$ meets the fibres of at most two atoms. Let $a$ be one of them, and let $y_a=\min(W\cap\rho^{-1}(a))$. If $\beta(y_a)\subseteq t_0$, then $y_a\vee p=a\in W$; so there is $i\in\beta(y_a)\setminus t_0\subseteq V_a$. Then $i$ lies in $\beta(x)$ for every $x\in W\cap\rho^{-1}(a)$, and for every $x\in W\cap L_2$, since $y_a\wedge x\in W\cap\rho^{-1}(a)$ gives $y_a\wedge x=y_a$ by minimality. If $W$ meets the fibre of a second atom $b$, we choose in the same way $j\in V_b$, contained in $\beta(x)$ for all $x\in W\cap(\rho^{-1}(b)\cup L_2)$. Otherwise let $w_0=\min(W\cap L_2)$. If $\beta(w_0)\subseteq\ut_a$, lifting $\beta(w_0)$ to $z$ in the fibre of $a$ gives $z\wedge b=p$ and $z\vee b=w_0$ for $b\in A_\dm\setminus\{a\}$, so $b\in W$; hence there is $j\in\beta(w_0)\setminus(t_0\cup V_a)$, contained in $\beta(x)$ for all $x\in W\cap L_2$.

In both cases $i,j\notin t_0$, $i\in V_a$ and $j\notin V_a$, so at most one of $i,j$ lies in $\beta(x)$ for $x\in L_1$. Hence $g=z_\dm-z_i-z_j$ is nonnegative on $L$, and it vanishes on $W$, where its values are $0$, $1-1$ and $2-1-1$ on the three levels. Its zero set is closed and contains no atom of the diamond, so it is distributive and equals $W$ by maximality. On $W$ we have $v_\dm=v_i+v_j$, and the only element of $[p,q]$ in $W$ is $p$; so $\beta$ is injective on $W$ by Lemma~\ref{lem:levels}(c), and Lemma~\ref{lem:boolean} for $\beta(W)\cong W$ shows that the vectors $(1,\one_{\beta(w)})$, and hence the vectors $(1,v(w))=(1,v_i(w)+v_j(w),\one_{\beta(w)})$, span a space of dimension $1+\rk W=1+\kappa(W)$. Proposition~\ref{prop:inclusion} gives $P_L\subseteq P_W$. If $q\in W$, apply this to $L^{\mathrm{op}}$, which has the same closed sets and the same ideals.
\end{proof}

\begin{theorem}\label{thm:atomfree}
Let $W$ be maximal among the closed distributive subsets of $L$, with $W\cap A_\dm=\emptyset$. Then $P_L\subseteq P_W$, unless $W\cap[p,q]=\emptyset$ and $\rk W=r-1$, in which case $P_W\in\Min(I_L)$ and $\htop P_W=n-r$.
\end{theorem}

\begin{proof}
The first case is Propositions~\ref{prop:endpoint} and~\ref{prop:avoid}. If $\rk W=r-1$, then $\htop P_W=n-1-\kappa(W)=n-r$, and $P_W$ is minimal because every prime containing $I_L$ has height at least $n-r$ (Proposition~\ref{prop:general}).
\end{proof}

\section{The minimal primes}\label{sec:main}

Recall the sets $\cD(L)$ and $\cN(L)$ of Definition~\ref{def:DN}. By Lemma~\ref{lem:M3}(b), each $D\in\cD(L)$ contains exactly one atom of the diamond, and each atom lies in some $D\in\cD(L)$, since $\{a\}$ is closed and distributive.

\begin{proof}[Proof of Theorem A]
The prime $P_L$ is minimal of height $n-1-\kappa(L)=n-r$ by Proposition~\ref{prop:minprimes}(c) and Lemma~\ref{lem:onediamond}. If $P_W\ne P_L$ is minimal, then $W$ does not contain the diamond by Proposition~\ref{prop:diamondinclusion}, so $W$ is distributive by Lemma~\ref{lem:M3}(b), and maximal by Lemma~\ref{lem:hibi}. Either $W$ meets $A_\dm$, and $W\in\cD(L)$, or $W\in\cN(L)$ by Theorem~\ref{thm:atomfree}. Conversely, the primes $P_W$ with $W\in\cN(L)$ are minimal by Theorem~\ref{thm:atomfree}. For $D\in\cD(L)$, let $P_W\subseteq P_D$ be a minimal prime. Then $D\subseteq W$ by Proposition~\ref{prop:minprimes}(a), $W$ contains only one atom by Lemma~\ref{lem:oneatom}, so it is distributive, and $W=D$ by maximality. The heights follow from Proposition~\ref{prop:minprimes}(c), since $\kappa(D)=\rk D\le r-1$, and distinct sets give distinct primes by Proposition~\ref{prop:minprimes}(a).
\end{proof}

In particular, all minimal primes of $I_L$ have the same height if and only if $\rk D=r-1$ for every $D\in\cD(L)$. Recall from Lemma~\ref{lem:rho} that the fibre of $c\in A_\dm$ under $\rho$ is $[m_c,M_c]$, with $\beta(m_c)=\lt_c$ and $\beta(M_c)=\ut_c$.

\begin{lemma}\label{lem:compl}
An atom $c\in A_\dm$ is complemented in its fibre if and only if $\lt_c\cup V_c\in\Delta$. In this case $s\cup(t\setminus t_0)\in\Delta$ for all $s,t\in\Delta$ with $\lt_c\subseteq s\subseteq t_0\subseteq t\subseteq\ut_c$.
\end{lemma}

\begin{proof}
A complement $c'$ of $c$ lies in $[m_c,M_c]$, and $\beta(c')\cap t_0=\lt_c$, $\beta(c')\cup t_0=\ut_c$, so $\beta(c')=\lt_c\cup V_c$. Conversely, let $c'$ be the element of the fibre of $c$ with $\beta(c')=\lt_c\cup V_c$ (Lemma~\ref{lem:rho}). Then $\beta(c'\wedge c)=\lt_c$ and $\beta(c'\vee c)=\ut_c$, and since $\beta$ is injective on the fibre, $c'\wedge c=m_c$ and $c'\vee c=M_c$. In this case $s\cup(t\setminus t_0)=(s\cup\beta(c'))\cap t\in\Delta$.
\end{proof}

\begin{proposition}\label{prop:local}
Suppose that some $c\in A_\dm$ is not complemented in its fibre. Then for every $b\in A_\dm\setminus\{c\}$ there is $D\in\cD(L)$ with $b\in D$ and $\rk D\le r-2$. In particular, the minimal primes of $I_L$ do not all have the same height.
\end{proposition}

\begin{proof}
Let $u=m_c\wedge p$ and $w=M_c\vee q$, so that $(\rho,\beta)(u)=(p,\lt_c)$ and $(\rho,\beta)(w)=(q,\ut_c)$. By Lemma~\ref{lem:compl}, $\lt_c\cup V_c\notin\Delta$; so $\lt_c\ne t_0$ and $V_c\ne\emptyset$, that is, $u<p$ and $q<w$. An element $x\in[u,w]$ in the fibre of an atom $e\ne c$ satisfies $\lt_c\cup\lt_e\subseteq\beta(x)\subseteq\ut_c\cap\ut_e$, so $\beta(x)=t_0$ and $x=e$. Hence $[u,b)\subseteq L_0$ and $(b,w]\subseteq L_2$.

We show that the chain $\{u,b,w\}$ is closed. Two incomparable elements with meet $u$ and join $b$ lie in $[u,b)\subseteq L_0$, so their join is not $b$; dually for meet $b$ and join $w$. Let $x,y$ be incomparable with $x\wedge y=u$ and $x\vee y=w$. Then $\rho(x)$ and $\rho(y)$ are distinct atoms, as $L_0$ lies below $L_2$. If neither is $c$, then $x$ and $y$ are atoms and $x\wedge y=p\ne u$. If $\rho(x)=c$, then $y$ is an atom, and $x\wedge y=u$, $x\vee y=w$ give $\beta(x)\cap t_0=\lt_c$ and $\beta(x)\cup t_0=\ut_c$, so $\beta(x)=\lt_c\cup V_c\notin\Delta$, a contradiction.

Let $D\supseteq\{u,b,w\}$ be maximal among the closed distributive sets, so that $D\in\cD(L)$, and let $e\in A_\dm\setminus\{b,c\}$. Since $M_c\wedge e=p$, $M_c\vee e=w$, $m_c\wedge e=u$ and $m_c\vee e=q$, neither $p$ nor $q$ lies in $D$, for otherwise $e\in D$. In a maximal chain of $D$ through $u,b,w$, the predecessor of $b$ lies strictly below $p$, and the successor strictly above $q$. So these two steps have length at least two in $L$, and $\rk D\le r-2$. By Theorem A, $\htop P_D\ge n-r+1>\htop P_L$.
\end{proof}

\begin{corollary}\label{cor:local}
Let $a\in A_\dm$, and let $y\lessdot a\lessdot z$ with $y\ne p$, $z\ne q$ and $[y,z]=\{y,a,z\}$. Then $a$ is not complemented in its fibre.
\end{corollary}

\begin{proof}
Since $y\lessdot a$ and $y\ne p$, $y\notin L_0$, so $\rho(y)=a$, and dually $\rho(z)=a$. If $a'$ were a complement of $a$ in its fibre, then $x=(y\vee a')\wedge z\in[y,z]$ would satisfy $x\wedge a=y\vee(a'\wedge a)=y$ and $x\vee a=(y\vee a'\vee a)\wedge z=z$ by the modular law, so $x\notin\{y,a,z\}$, a contradiction.
\end{proof}

\section{Radicality}\label{sec:radical}

Grade $S_L$ by $\deg x_a=\Phi_L(e_a)=(1,v(a))$. For a degree $m$, the set $F_m=\{\alpha\in\N^L:\Phi_L(\alpha)=m\}$ is finite, and we identify $x^\alpha$ with the multiset $\alpha$. A \emph{move} replaces an incomparable pair by its meet and join, or conversely. We write $\alpha\sim\beta$ if $\alpha$ and $\beta$ are connected by moves, and $\mathcal M_m=F_m/{\sim}$.

\begin{lemma}\label{lem:moves}
Let $L$ be a finite modular lattice and $K$ any field.
\begin{enumerate}
\item[(a)] $I_L$ and $P_W$, for every closed $W$, are homogeneous.
\item[(b)] $x^\alpha-x^\beta\in I_L$ if and only if $\alpha\sim\beta$, and the classes in $\mathcal M_m$ give a basis of $(S_L/I_L)_m$.
\item[(c)] Let $W$ be closed and $\alpha\sim\beta$. Then $\alpha$ is supported in $W$ if and only if $\beta$ is. If $\alpha$ is not supported in $W$, then $x^\alpha\in P_W$. Two monomials supported in $W$ have the same image in $S_L/P_W$ if and only if they have the same image under $\Phi_W$, and this image is nonzero.
\end{enumerate}
\end{lemma}

\begin{proof}
(a) The join-meet generators are homogeneous by Lemma~\ref{lem:valuation}, and so are the generators of $Q_W$, since $\Lambda_W\subseteq\ker\Phi_L$. (b) In degree $m$, the products of generators with monomials span the differences of monomials joined by a move, that is, the kernel of $K^{F_m}\to K^{\mathcal M_m}$ sending each monomial to its class; and a binomial of two different degrees is not in $I_L$, since $I_L$ is homogeneous and contains no monomial. (c) Moves preserve support in $W$ because $W$ is closed, and the rest follows from $R_W/Q_W\cong\operatorname{im}\phi_W$.
\end{proof}

\begin{theorem}\label{thm:radcrit}
Let $L$ be a finite modular lattice with exactly one diamond and $K$ any field. For a degree $m$, let $\mathcal R_m\subseteq\mathcal M_m$ be the set of classes none of whose elements is supported on a closed distributive subset of $L$, and $e_m=|\mathcal R_m|$. Then $\dim_K(\sqrt{I_L}/I_L)_m=\max\{e_m-1,0\}$. In particular, $I_L$ is radical if and only if $e_m\le1$ for every $m$, and this does not depend on $K$.
\end{theorem}

\begin{proof}
By Proposition~\ref{prop:minprimes} and Lemma~\ref{lem:moves}, $(\sqrt{I_L}/I_L)_m$ is the kernel of $K^{\mathcal M_m}\to\bigoplus_W(S_L/P_W)_m$, $W$ closed. If a class $C$ contains some $\alpha$ supported on a closed distributive $W$, no other class has an element supported on $W$ with the same $\Phi_W$-degree, since Hibi's theorem would give $x^\alpha-x^\beta\in I_W\subseteq I_L$; so, by Lemma~\ref{lem:moves}(c), the image of $x^\alpha$ in $S_L/P_W$ shows that the coefficient of $C$ in a vector of the kernel vanishes. The classes in $\mathcal R_m$ vanish modulo every $P_W$ with $W$ distributive. If $W$ is not distributive, then $P_L\subseteq P_W$ by Proposition~\ref{prop:diamondinclusion}; all monomials of degree $m$ have the same nonzero image modulo $P_L$, so each such $W$ imposes at most the condition that the coefficients add up to zero, and $W=L$ imposes it. Hence the kernel consists of the vectors supported on $\mathcal R_m$ with coefficient sum zero.
\end{proof}

We now prove Theorem B: sufficiency in Section~\ref{subsec:suff}, necessity in Section~\ref{subsec:nec}. Corollary C is proved at the end of the section.

\subsection{Defects}
From now on $L$ has exactly one diamond. The \emph{defect} of $x\in L$ is $\delta(x)=\beta(x)\mathbin{\triangle}t_0$, and we put $\delta^-(x)=\delta(x)\cap t_0$ and $\delta^+(x)=\delta(x)\setminus t_0$. By Lemma~\ref{lem:rho}, $x$ is determined by $\rho(x)$ and $\delta(x)$, and we write $x=\langle\rho(x),\delta(x)\rangle$. By Lemma~\ref{lem:beta}(a), $\delta^-$ turns meets into unions and joins into intersections, and $\delta^+$ does the opposite. Elements of $L_0$ have $\delta=\delta^-$, elements of $L_2$ have $\delta=\delta^+$, and elements of $\rho^{-1}(c)$ have $\delta^-\subseteq R_c$ and $\delta^+\subseteq V_c$.

We say that $E\in\cE$ has \emph{color} $c$ if $E\in R_c\cup V_c$. Write $\col(Z)$ for the set of colors occurring in $Z\subseteq\cE$, and call $Z$ \emph{bicolored} if all its elements have a color and $|\col(Z)|\le2$. For a multiset of degree $m$ and size $N$, let $\nu_E=\nu_E(m)$ count the occurrences of $E$ in its defects, and $d=d(m)=\sum_xv_\dm(x)$. Since $\nu_E=N-\sum_xv_E(x)$ for $E\in t_0$ and $\nu_E=\sum_xv_E(x)$ otherwise, the numbers $N$, $d$ and $\nu_E$ determine $m$, and conversely.

\begin{lemma}\label{lem:split}
Let $u\in L_0$, $w\in L_2$ and $Z=\delta(u)\cup\delta(w)$.
\begin{enumerate}
\item[(a)] If $x,y$ are incomparable with $x\wedge y=u$ and $x\vee y=w$, then $\rho(x)$ and $\rho(y)$ are distinct atoms $e$ and $f$, $\delta(x)=Z\cap(R_e\cup V_e)$ and $\delta(y)=Z\cap(R_f\cup V_f)$. In particular $Z$ is bicolored.
\item[(b)] Suppose that every atom of the diamond is complemented in its fibre and that $Z$ is bicolored. Let $e,f\in A_\dm$ be distinct with $\col(Z)\subseteq\{e,f\}$. Then $x_e=\langle e,Z\cap(R_e\cup V_e)\rangle$ and $x_f=\langle f,Z\cap(R_f\cup V_f)\rangle$ exist, they are incomparable, $x_e\wedge x_f=u$ and $x_e\vee x_f=w$.
\end{enumerate}
\end{lemma}

\begin{proof}
(a) The images under $\rho$ have meet $p$ and join $q$, and they are distinct atoms, since $L_0$ lies below $L_2$. The defects of $x$ and $y$ have disjoint colors, the meet takes the union of their negative parts and the join the union of their positive parts. (b) For $e$, put $s=\beta(u)\cup\lt_e$ and $t=\beta(w)\cap\ut_e$. By Lemma~\ref{lem:compl}, $s\cup(t\setminus t_0)\in[\lt_e,\ut_e]_\Delta$, and its lift to the fibre of $e$ is $x_e$; similarly for $x_f$. The rules for defects give the meet and the join.
\end{proof}

\subsection{Multichains}
Fix a degree $m$ with total degree $N\ge1$, and write $\nu_E=\nu_E(m)$ and $d=d(m)$. For $1\le i\le N$ put
\[
\delta_i=\{E\in t_0 : \nu_E\ge i\}\cup\{E\in\cE\setminus t_0 : \nu_E\ge N-i+1\}.
\]
For an integer $k$ and $c\in A_\dm$ put $l=d-2k$ and $n_0=N-l-k$. If $n_0,l,k\ge0$ and the elements $\langle p,\delta_i\rangle$ for $1\le i\le n_0$, $\langle c,\delta_i\rangle$ for $n_0<i\le n_0+l$ and $\langle q,\delta_i\rangle$ for $n_0+l<i\le N$ exist, we denote by $\zeta_m(k,c)$ the multiset formed by them; it does not depend on $c$ if $l=0$. Each $E$ lies in exactly $\nu_E$ of the sets $\delta_i$, so $\zeta_m(k,c)$ has degree $m$, and it is a multichain, because $\rho$ and $\delta^+$ increase and $\delta^-$ decreases along it.

\begin{lemma}\label{lem:chains}
\begin{enumerate}
\item[(a)] Every multiset is connected by moves to a multichain, and every multichain of degree $m$ is of the form $\zeta_m(k,c)$.
\item[(b)] If $\zeta_m(k,c)$ and $\zeta_m(k',c')$ exist and $k\le j\le k'$, then $\zeta_m(j,c)$ exists.
\item[(c)] If the class of $\zeta_m(k,c)$ lies in $\mathcal R_m$, then $n_0\ge1$ and $k\ge1$.
\end{enumerate}
\end{lemma}

\begin{proof}
(a) For incomparable $x,y$ put $a=\rk x-\rk(x\wedge y)$ and $b=\rk y-\rk(x\wedge y)$. Since $\rk(x\vee y)=\rk(x\wedge y)+a+b$, replacing $x,y$ by $x\wedge y,x\vee y$ increases the sum of the squared ranks by $2ab>0$; as this sum is bounded, the process ends at a multichain. Along a multichain the negative defects decrease and the positive defects increase, which gives the sets $\delta_i$, and the $\rho$-values are $p$, then one atom $c$, then $q$; so it is $\zeta_m(k,c)$, where $k$ is the number of elements in $L_2$.

(b) Put $n_0(j)=N-d+j$. For $i\le n_0(j)\le n_0(k')$ the element $\langle p,\delta_i\rangle$ occurs in $\zeta_m(k',c')$, and so does $\langle q,\delta_i\rangle$ for $i>N-j\ge N-k'$; for the remaining indices $n_0(k)<i\le N-k$, so $\langle c,\delta_i\rangle$ occurs in $\zeta_m(k,c)$. (c) If $k=0$, the multichain lies in the closed interval ${\downarrow}z_N$, which does not contain $q$ and is therefore distributive; the case $n_0=0$ is dual.
\end{proof}

\begin{lemma}\label{lem:branch}
Let $c\in A_\dm$, and let two multisets of the same degree be supported in $\rho^{-1}(\{p,c,q\})$ and have the same number of elements in $L_2$. Then they are connected by moves.
\end{lemma}

\begin{proof}
Since $\{p,c,q\}$ is a chain, a move inside the sublattice $\rho^{-1}(\{p,c,q\})$ preserves the numbers of elements in $L_0$, in the fibre of $c$ and in $L_2$. As in Lemma~\ref{lem:chains}(a), both multisets are connected inside $\rho^{-1}(\{p,c,q\})$ to multichains $z_1\le\cdots\le z_N$ of the same degree with the same number $k$ of elements in $L_2$. In such a multichain $\rho(z_i)$ is $p$, $c$ or $q$ according as $i\le N-d+k$, $N-d+k<i\le N-k$ or $i>N-k$, and for $E\in\cE$ the indices $i$ with $E\in\beta(z_i)$ form a final segment of length $\sum_iv_E(z_i)$. So the multichain is determined by the degree and by $k$, since $(\rho,\beta)$ is injective.
\end{proof}

\subsection{Sufficiency}\label{subsec:suff}

\begin{lemma}\label{lem:W}
Let $c\in A_\dm$, $X\subseteq t_0\setminus R_c$ and $Y\subseteq\cE\setminus(t_0\cup V_c)$, and suppose that $X\cup Y$ is not bicolored. Then
\[
W=\{x\in L_0 : \delta(x)\supseteq X\}\cup\rho^{-1}(c)\cup\{x\in L_2 : \delta(x)\supseteq Y\}
\]
is a closed distributive subset of $L$.
\end{lemma}

\begin{proof}
The three parts of $W$ are sublattices, by the rules for defects, and if $x\in W\cap L_0$ and $\rho(y)=c$, then $x\wedge y\in L_0$ has defect $\delta(x)\cup\delta^-(y)\supseteq X$ and $x\vee y\in\rho^{-1}(c)$; with the dual remark for $L_2$, and since elements of $L_0$ and $L_2$ are comparable, $W$ is a sublattice. Let $x,y$ be incomparable with $x\wedge y,x\vee y\in W$. If $\rho(x)$ is an atom $e\ne c$, then $\rho(y)$ is another atom, $x\wedge y\in L_0$, $x\vee y\in L_2$, and Lemma~\ref{lem:split}(a) shows that $\delta(x\wedge y)\cup\delta(x\vee y)\supseteq X\cup Y$ is bicolored, a contradiction. So $\rho(x),\rho(y)\in\{p,c,q\}$, and we may assume $x\notin\rho^{-1}(c)\subseteq W$. Let $x\in L_0$, so that $y\notin L_2$. If $y\in L_0$, then $\delta(x)\supseteq\delta(x\vee y)\supseteq X$; if $\rho(y)=c$, then $X\subseteq\delta(x)\cup\delta^-(y)$ with $\delta^-(y)\subseteq R_c$ disjoint from $X$. In both cases $x\in W$, and the case $x\in L_2$ is dual. So $W$ is closed, and it is distributive by Lemma~\ref{lem:M3}(b), since $W\cap A_\dm=\{c\}$.
\end{proof}

In the following three lemmas we assume that every atom of the diamond is complemented in its fibre.

\begin{lemma}\label{lem:transfer}
Let $\zeta=\zeta_m(k,c)=\{z_1\le\dots\le z_N\}$ be a multichain whose class lies in $\mathcal R_m$, with $l=d-2k\ge1$ and $n_0=N-l-k$. Suppose that $\delta_{n_0+1}\subseteq t_0$ and $\delta_{n_0+l}\cap t_0=\emptyset$, and put $X=\delta_{n_0}\setminus R_c$ and $Y=\delta_{n_0+l+1}\setminus V_c$. Then:
\begin{enumerate}
\item[(a)] $X\cup Y$ is bicolored, and $c\notin\col(X\cup Y)$;
\item[(b)] the elements $u=\langle p,X\rangle$, $w=\langle q,Y\rangle$ and, for $1\le j\le l$,
\[
y_j=\bigl\langle c,\ (\delta_{n_0+j-1}\cap R_c)\cup(\delta_{n_0+j+1}\cap V_c)\bigr\rangle
\]
exist, and $\zeta\sim\Theta$, where
\[
\Theta=\{z_1,\dots,z_{n_0-1},u,y_1,\dots,y_l,w,z_{n_0+l+2},\dots,z_N\}.
\]
\end{enumerate}
\end{lemma}

\begin{proof}
Lemma~\ref{lem:chains}(c) gives $n_0,k\ge1$. (a) Otherwise Lemma~\ref{lem:W} gives a closed distributive set supporting $\zeta$, since the defects of the elements of $\zeta$ in $L_0$ contain $\delta_{n_0}\supseteq X$, those in $L_2$ contain $\delta_{n_0+l+1}\supseteq Y$, and the other elements lie in $\rho^{-1}(c)$. The color $c$ does not occur, since $X\cap R_c=Y\cap V_c=\emptyset$.

(b) By (a), all elements of $\delta_{n_0}$ have a color, and $\beta(u)=t_0\setminus X=\bigcap_{e\ne c}(\beta(z_{n_0})\cup\lt_e)$ and $\beta(w)=t_0\cup Y=\bigcup_{e\ne c}(\beta(z_{n_0+l+1})\cap\ut_e)$ lie in $\Delta$; we get $u$ and $w$ by meeting a preimage of the first with $p$ and joining a preimage of the second with $q$. For $y_j$, apply Lemma~\ref{lem:compl} to $s=(\beta(z_{n_0+j-1})\cap t_0)\cup\lt_c$ and $t=(\beta(z_{n_0+j+1})\cup t_0)\cap\ut_c$. The multisets $\zeta$ and $\Theta$ have $n_0$ elements in $L_0$, $l$ in the fibre of $c$ and $k$ in $L_2$, so by Lemma~\ref{lem:branch} it suffices to show that each $E\in\cE$ occurs equally often in the defects of their factors with indices $n_0,\dots,n_0+l+1$. Labels outside $R_c\cup V_c$ occur there only through $\delta_{n_0}\setminus R_c=X=\delta(u)$ and $\delta_{n_0+l+1}\setminus V_c=Y=\delta(w)$. For the labels in $R_c$, passing from $\zeta$ to $\Theta$ shifts the sets $\delta_{n_0}\cap R_c,\dots,\delta_{n_0+l}\cap R_c$ one place to the right, and the last one, which disappears, is empty because $\delta_{n_0+l}\cap t_0=\emptyset$. Dually, the sets $\delta_{n_0+1}\cap V_c,\dots,\delta_{n_0+l+1}\cap V_c$ are shifted one place to the left, and $\delta_{n_0+1}\cap V_c=\emptyset$.
\end{proof}

\begin{lemma}\label{lem:raise}
Let $\zeta_m(k,c)$ be a multichain whose class lies in $\mathcal R_m$, and suppose that $d-2k\ge2$ and that $\zeta_m(k+1,c)$ exists. Then $\zeta_m(k,c)\sim\zeta_m(k+1,c)$.
\end{lemma}

\begin{proof}
Put $l=d-2k$ and $n_0=N-l-k$. In $\zeta_m(k+1,c)$ the elements with indices $n_0+1$ and $n_0+l$ lie in $L_0$ and $L_2$, so $\delta_{n_0+1}\subseteq t_0$, $\delta_{n_0+l}\cap t_0=\emptyset$, and Lemma~\ref{lem:transfer} applies. Since $|A_\dm|\ge3$, we can choose distinct $e,f\in A_\dm\setminus\{c\}$ with $\col(X\cup Y)\subseteq\{e,f\}$, and replace $u,w$ in $\Theta$ by $x_e,x_f$ (Lemma~\ref{lem:split}(b)). The pairs $x_e,y_1$ and $x_f,y_2$ are incomparable, since their images under $\rho$ are; replacing them by their meets and joins we get a multiset supported in $\rho^{-1}(\{p,c,q\})$ with $k+1$ elements in $L_2$, which is connected to $\zeta_m(k+1,c)$ by Lemma~\ref{lem:branch}.
\end{proof}

\begin{lemma}\label{lem:top}
Suppose that $\zeta_m(k,c)$ has exactly one element in $L_1$, and that this element is $c$. Then $\zeta_m(k,c')$ exists for every $c'\in A_\dm$, and all the multichains $\zeta_m(k,c')$ whose classes lie in $\mathcal R_m$ lie in the same class.
\end{lemma}

\begin{proof}
The middle defect is empty, so we may replace $c$ by any $c'$. Let $\mathcal A$ be the set of $c'\in A_\dm$ such that the class of $\zeta_m(k,c')$ lies in $\mathcal R_m$, and assume $\mathcal A\ne\emptyset$; then $n_0=N-k-1\ge1$ and $k\ge1$ by Lemma~\ref{lem:chains}(c). Put $Z=\delta_{n_0}\cup\delta_{n_0+2}$ and $C=\col(Z)$. For $a\in\mathcal A$, Lemma~\ref{lem:transfer} says that $Z\setminus(R_a\cup V_a)$ is bicolored; so every element of $Z$ has a color and $|C\setminus\{a\}|\le2$.

We claim that if $e,f\in A_\dm\setminus\{a\}$ are distinct and $C\setminus\{a\}\subseteq\{e,f\}$, then $\zeta_m(k,a)\sim\zeta_m(k,e)\sim\zeta_m(k,f)$. Indeed, take $\Theta$ as in Lemma~\ref{lem:transfer} for $\zeta_m(k,a)$, replace $u,w$ by $x_e,x_f$, and then $y_1,x_f$ by their meet and join. The result is supported in $\rho^{-1}(\{p,e,q\})$ and has $k$ elements in $L_2$, so Lemma~\ref{lem:branch} connects it to $\zeta_m(k,e)$; exchanging $e$ and $f$ gives the other connection. In particular $e,f\in\mathcal A$.

If $|C\setminus\{a\}|\le1$ for some $a\in\mathcal A$, let $g\in A_\dm\setminus\{a\}$, and take $f\in C\setminus\{a\}$ if this set is nonempty and different from $\{g\}$, and any $f\in A_\dm\setminus\{a,g\}$ otherwise; the claim with $e=g$ gives $\zeta_m(k,a)\sim\zeta_m(k,g)$. Otherwise $|C\setminus\{a\}|=2$ for all $a\in\mathcal A$. Fix $a$ and write $C\setminus\{a\}=\{e,f\}$; then $e,f\in\mathcal A$ by the claim. If $a\notin C$, then $C=\{e,f\}$ and $|C\setminus\{e\}|=1$, a contradiction. So $C=\{a,e,f\}$, every element of $\mathcal A$ lies in $C$, since otherwise removing it leaves three colors, and the claim connects the three multichains.
\end{proof}

\begin{proof}[Proof of Theorem B, sufficiency]
By Theorem~\ref{thm:radcrit} we must show $e_m\le1$. Two classes in $\mathcal R_m$ contain multichains $\zeta_m(k,c)$ and $\zeta_m(k',c')$ by Lemma~\ref{lem:chains}(a). Let $k^*$ be the largest $j$ for which some $\zeta_m(j,\cdot)$ exists. By Lemma~\ref{lem:chains}(b), $\zeta_m(j,c)$ exists for $k\le j\le k^*$, and $d-2j\ge2$ for $j<k^*$, so Lemma~\ref{lem:raise} gives $\zeta_m(k,c)\sim\zeta_m(k^*,c)$; similarly $\zeta_m(k',c')\sim\zeta_m(k^*,c')$. If these two multichains differ, then $l=d-2k^*\ge1$, $c\ne c'$, and the defects of their elements in $L_1$ lie in $(R_c\cup V_c)\cap(R_{c'}\cup V_{c'})=\emptyset$. If $l\ge2$, replacing the first and the last of these atoms by $p$ and $q$ gives $\zeta_m(k^*+1,c)$, against the choice of $k^*$. So $l=1$, and Lemma~\ref{lem:top} shows that the two classes coincide.
\end{proof}

\subsection{Necessity}\label{subsec:nec}

\begin{lemma}\label{lem:ij}
Let $c\in A_\dm$ be an atom that is not complemented in its fibre. Then there are $i\in R_c$ and $j\in V_c$ such that every element of $\Delta$ containing $j$ contains $i$.
\end{lemma}

\begin{proof}
For $Z\subseteq\cE$ let $\overline Z=\bigcap\{t\in\Delta:Z\subseteq t\}\in\Delta$. By Lemma~\ref{lem:compl}, $h=\lt_c\cup V_c\notin\Delta$, so $\overline h\ne h$; since $\overline h\subseteq\ut_c$ and $\ut_c\setminus h=R_c$, there is $i\in\overline h\cap R_c$. If $i\notin\overline{\{j\}}$ for every $j\in V_c$, then $\lt_c\cup\bigcup_{j\in V_c}\overline{\{j\}}$ is an element of $\Delta$ containing $h$ but not $i$, a contradiction.
\end{proof}

\begin{proposition}\label{prop:witness}
Let $c\in A_\dm$ be an atom that is not complemented in its fibre, and put $u=m_c\wedge p$ and $w=M_c\vee q$. Then $\alpha=\{u,c,M_c,w\}$ and $\alpha'=\{u,p,w,w\}$ have the same degree $m$, their classes lie in $\mathcal R_m$, and $\alpha\not\sim\alpha'$. Moreover, over every field, $F=x_ux_w(x_cx_{M_c}-x_px_w)\in\sqrt{I_L}\setminus I_L$.
\end{proposition}

\begin{proof}
By Lemma~\ref{lem:compl}, $\lt_c\ne t_0$ and $V_c\ne\emptyset$, so $m_c\ne c\ne M_c$. Let $e\in A_\dm\setminus\{c\}$. By Lemma~\ref{lem:rho}, the pairs $M_c,e$, then $M_c,q$, then $m_c,p$ are incomparable, with
\[
M_c\wedge e=p,\ M_c\vee e=w;\qquad M_c\wedge q=c,\ M_c\vee q=w;\qquad m_c\wedge p=u,\ m_c\vee p=c.
\]
The first pair and $v(e)=v(c)$ show that $\alpha$ and $\alpha'$ have the same degree. A closed set containing $\alpha$ contains $c,w$, hence $q$, and $u,c$, hence $p$; one containing $\alpha'$ contains $p,w$, hence every $e\ne c$. By Lemmas~\ref{lem:M3}(b) and~\ref{lem:moves}(c), both classes lie in $\mathcal R_m$.

In $S_L/I_L$ we write $t$ for the class of $x_t$, put $y=m_c$, $z=M_c$, and choose distinct $a,b\in A_\dm\setminus\{c\}$. The relations $za=zb=pw$, $zq=cw$, $uc=yp$ and $ca=cb=ab=pq$ give, for $H=cz-pw=z(c-a)=z(c-b)$,
\[
H^2=z^2(c-a)(c-b)=z^2(c^2-pq)=czH,\qquad bH=z(bc-ab)=0,
\]
where $z^2pq=z(zq)p=czpw$. So the class of $F$ is $uwH$, and its square is $uwcz\cdot uwH=ypwz\cdot uwH=yz^2b\cdot uwH=0$.

It remains to show $\alpha\not\sim\alpha'$, so that $F\notin I_L$ by Lemma~\ref{lem:moves}(b). Choose $i,j$ as in Lemma~\ref{lem:ij}, and let $\psi:L\to[p,q]\times2^{\{i,j\}}$, $\psi(x)=(\rho(x),\beta(x)\cap\{i,j\})$. A move of $L$ is sent by $\psi$ to a move of the lattice $\psi(L)$ or to the identity, so it suffices to show that $\psi(\alpha)$ and $\psi(\alpha')$ are not connected in $\psi(L)$. We write $00,10,11$ for $\emptyset,\{i\},\{i,j\}$. Elements of $L_0$ have pattern $00$ or $10$, elements of $L_2$ have $10$ or $11$, the fibre of $c$ has the patterns $00,10,11$ by the choice of $i,j$, and the fibre of $e\ne c$ has pattern $10$, since $i\in R_c\subseteq\lt_e$ and $j\notin\ut_e$. Put $P_0=(p,00)$, $P_1=(p,10)$, $C_1=(c,10)$, $C_2=(c,11)$, $Q_2=(q,11)$ and $E_e=(e,10)$; then $\psi(\alpha')=P_0P_1Q_2^2$ and $\psi(\alpha)=P_0C_1C_2Q_2$. Two incomparable elements of $\psi(L)$ with join $Q_2$ and meet of level $0$ cannot have $\rho$-values $p$ and $q$, since every pattern at level $0$ is contained in every pattern at level $2$; so they lie in the fibres of two distinct atoms, one of them is $C_2$, as it must contain $j$, the other is some $E_e$, and their meet is $P_1$. Two incomparable elements with meet $P_0$ and join $C_2$ do not exist: they would have $\rho$-values in $\{p,c\}$, the patterns in the fibre of $c$ form a chain, and if one of them has level $0$, the other contains $j$ and is $C_2$. Since moreover the intervals $[P_0,P_1]$, $[P_0,E_e]$, $[C_2,Q_2]$ and $[E_e,Q_2]$ of $\psi(L)$ are chains, the only moves from $P_0P_1Q_2^2$ lead to the monomials $P_0C_2E_eQ_2$, and the only moves from these lead back. So $\psi(\alpha)$ is not connected to $\psi(\alpha')$.
\end{proof}

\begin{proof}[Proof of Theorem B, necessity]
This is Proposition~\ref{prop:witness}.
\end{proof}

\begin{proof}[Proof of Corollary C]
Minimal primes are associated primes, so (a) implies (b), and (b) is equivalent to (c) by Theorem A. Under (c), Proposition~\ref{prop:local} shows that every atom of the diamond is complemented in its fibre, so $I_L$ is radical by Theorem B; then its associated primes are its minimal primes, and (b) implies (a).
\end{proof}

\section{Examples and conjectures}\label{sec:final}

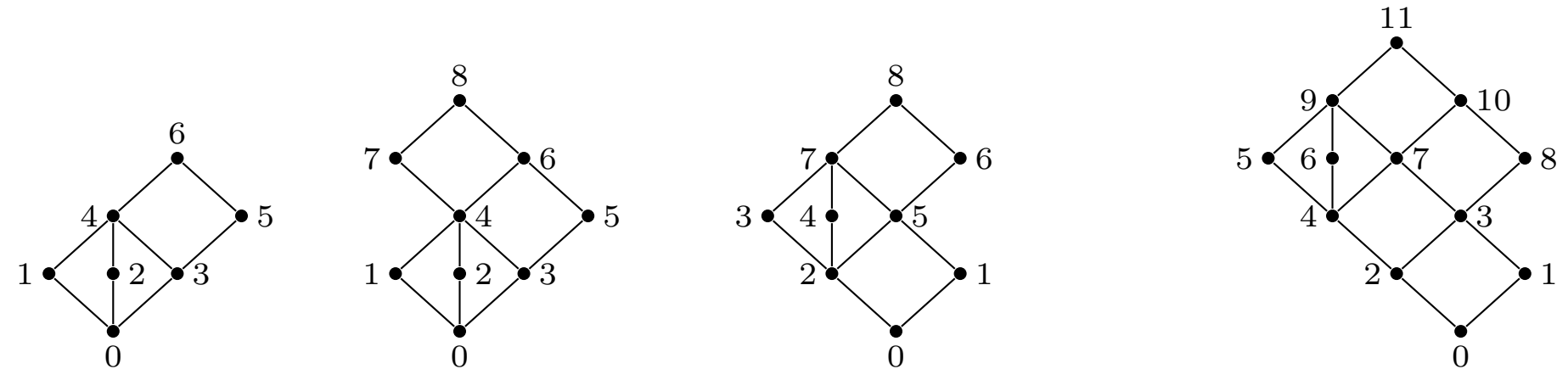
\begin{figure}[htbp]
\centering
\begin{tikzpicture}[x=0.5cm,y=0.45cm]
\begin{scope}
\foreach \n/\x/\y/\p in {0/0/0/below,1/-1/1/left,2/0/1/right,3/1/1/right,4/0/2/left,5/2/2/right,6/1/3/above} \node[dot,label={[lab]\p:$\n$}] (A\n) at (\x,\y) {};
\foreach \s/\t in {0/1,0/2,0/3,1/4,2/4,3/4,3/5,4/6,5/6} \draw (A\s)--(A\t);
\node[lab,align=center] at (0.5,-1.6) {(a) radical,\\ not unmixed};
\end{scope}
\begin{scope}[xshift=2.7cm]
\foreach \n/\x/\y/\p in {0/0/0/below,1/-1/1/left,2/0/1/right,3/1/1/right,4/0/2/right,5/2/2/right,6/1/3/right,7/-1/3/left,8/0/4/above} \node[dot,label={[lab]\p:$\n$}] (Z\n) at (\x,\y) {};
\foreach \s/\t in {0/1,0/2,0/3,1/4,2/4,3/4,3/5,4/6,5/6,4/7,6/8,7/8} \draw (Z\s)--(Z\t);
\node[lab,align=center] at (0.5,-1.6) {(b) unmixed\\ and radical};
\end{scope}
\begin{scope}[xshift=6.1cm]
\foreach \n/\x/\y/\p in {0/0/0/below,1/1/1/right,2/-1/1/left,3/-2/2/left,4/-1/2/left,5/0/2/right,6/1/3/right,7/-1/3/left,8/0/4/above} \node[dot,label={[lab]\p:$\n$}] (B\n) at (\x,\y) {};
\foreach \s/\t in {0/1,0/2,1/5,2/3,2/4,2/5,3/7,4/7,5/7,5/6,6/8,7/8} \draw (B\s)--(B\t);
\node[lab,align=center] at (-0.5,-1.6) {(c) not radical,\\ not unmixed};
\end{scope}
\begin{scope}[xshift=10.5cm]
\foreach \n/\x/\y/\p in {0/0/0/below,1/1/1/right,2/-1/1/left,3/0/2/right,4/-2/2/left,5/-3/3/left,6/-2/3/left,7/-1/3/right,8/1/3/right,9/-2/4/left,10/0/4/right,11/-1/5/above} \node[dot,label={[lab]\p:$\n$}] (L\n) at (\x,\y) {};
\foreach \s/\t in {0/1,0/2,1/3,2/3,2/4,3/7,3/8,4/5,4/6,4/7,5/9,6/9,7/9,7/10,8/10,9/11,10/11} \draw (L\s)--(L\t);
\node[lab,align=center] at (-1,-1.6) {(d) not radical,\\ not unmixed};
\end{scope}
\end{tikzpicture}
\caption{The lattices $A$, $Z$, $B$ and $L$ of Examples~\ref{ex:A}, \ref{ex:Z}, \ref{ex:B} and~\ref{ex:L12}, with the properties of their join-meet ideals.}\label{fig:examples}
\end{figure}

We illustrate our results with some examples and then propose some conjectures. In Examples~\ref{ex:A} and~\ref{ex:Z} the join-meet ideal is radical; they show that the converse of Proposition~\ref{prop:local} fails and that unmixedness is not inherited by intervals. In Examples~\ref{ex:B} and~\ref{ex:L12} the join-meet ideal is not radical. In the first one the atom without a complement in its fibre is detected by the local criterion of Corollary~\ref{cor:local}, while in the second one it is not, so this criterion is sufficient but not necessary. The lattices of Figure~\ref{fig:examples} have no sublattice isomorphic to $N_5$, so they are modular, and each of them has exactly one diamond.

\begin{example}\label{ex:A}
Let $A$ be the lattice of Figure~\ref{fig:examples}(a); here $n=7$ and $r=3$. The fibres of the atoms $1,2,3$ are $\{1\}$, $\{2\}$ and $\{3,5\}$, and $3$ is complemented by $5$ in its fibre, so $I_A$ is radical by Theorem B. The chain $\{0,1\}$ is closed. It is maximal among the closed distributive sets, because a closed set strictly containing it contains a second atom of $[0,4]$: if it contains $4$, it contains $2$ and $3$, as $2\wedge3=0$ and $2\vee3=4$; if it contains $5$, it contains $1\vee5=6$; and if it contains $6$, it contains $2$, as $2\wedge5=0$ and $2\vee5=6$. So $\{0,1\}\in\cD(A)$ has rank $1<r-1$, and by Theorem A the minimal primes of $I_A$ have heights $4$ and $5$. In particular $I_A$ is not unmixed, and the converse of Proposition~\ref{prop:local} fails.
\end{example}

\begin{example}\label{ex:Z}
The lattice $Z$ of Figure~\ref{fig:examples}(b) contains $A=[0,6]$ as an interval, and $r=4$. The chain $\{0,1\}$ is no longer maximal, since $\{0,1,7,8\}$ is a closed chain. A computer check shows that all elements of $\cD(Z)$ have rank $3$, so $I_Z$ is unmixed and radical by Corollary C; Macaulay2 \cite{M2} shows that it is even Cohen--Macaulay. So unmixedness of join-meet ideals is not inherited by intervals.
\end{example}

\begin{example}\label{ex:B}
In the lattice $B$ of Figure~\ref{fig:examples}(c), of rank $4$, the diamond is $[2,7]$ and $1\lessdot5\lessdot6$ is a configuration as in Corollary~\ref{cor:local}. So the atom $5$ is not complemented in its fibre $\{1,5,6\}$, and $I_B$ is neither radical nor unmixed. With $u=0$ and $w=8$ as in Proposition~\ref{prop:local}, the closed chain $\{0,3,8\}$ lies in $\cD(B)$ and has rank $2$, and the minimal primes of $I_B$ have heights $5$ and $6$. Proposition~\ref{prop:witness} gives $x_0x_8(x_5x_6-x_2x_8)\in\sqrt{I_B}\setminus I_B$.
\end{example}

\begin{example}\label{ex:L12}
In the lattice $L$ of Figure~\ref{fig:examples}(d), of rank $5$, the diamond is $[4,9]$, and the configuration of Corollary~\ref{cor:local} does not occur: the atoms $5$ and $6$ only cover $4$ and are only covered by $9$, and the only choice for $7$ is $3\lessdot7\lessdot10$, but $[3,10]=\{3,7,8,10\}$. Nevertheless the fibre of $7$ is $\{1,3,7,8,10\}$, and $7$ has no complement in it, since $7\vee1=7$ and $7\wedge x\in\{3,7\}$ for $x\in\{3,8,10\}$. So $I_L$ is neither radical nor unmixed, and $x_0x_{11}(x_7x_{10}-x_4x_{11})\in\sqrt{I_L}\setminus I_L$. The closed chain $D=\{0,5,11\}$ lies in $\cD(L)$ and has rank $2=r-3$: a closed distributive set containing $D$ cannot contain $6$ or $7$, the pairs $6,8$, $6,10$ and $1,6$ exclude $2$, $4$ and $9$, and then $1,3,8,10$ are excluded because $1\vee5=9$, $3\wedge5=8\wedge5=2$ and $10\wedge5=4$. The minimal primes of $I_L$ have heights $7$, $8$ and $9$.
\end{example}

In the next example we fix a planar distributive lattice and create the diamond by adding an atom to one of its squares; moving the diamond along the lattice, we obtain all the behaviours of $I_L$ allowed by Theorem B and Corollary C.

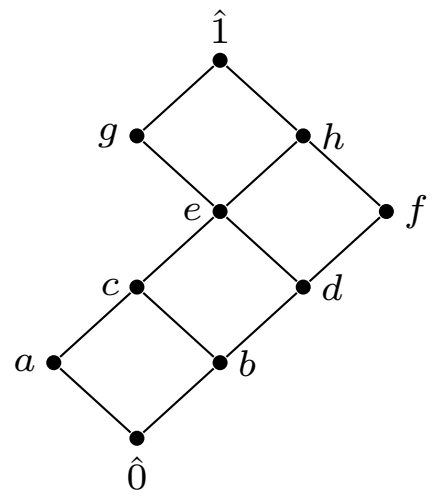
\begin{figure}[htbp]
	\centering
	\begin{tikzpicture}[x=0.55cm,y=0.5cm]
		\foreach \n/\l/\x/\y/\p in {o/\hat0/0/0/below,a/a/-1/1/left,b/b/1/1/right,c/c/0/2/left,d/d/2/2/right,e/e/1/3/left,f/f/3/3/right,g/g/0/4/left,h/h/2/4/right,i/\hat1/1/5/above} \node[dot,label={[lab]\p:$\l$}] (H\n) at (\x,\y) {};
		\foreach \s/\t in {o/a,o/b,a/c,b/c,b/d,c/e,d/e,d/f,e/g,e/h,f/h,g/i,h/i} \draw (H\s)--(H\t);
	\end{tikzpicture}
	\caption{The planar distributive lattice $L_0$ of Example~\ref{ex:hibi}.}\label{fig:hibi}
\end{figure}

\begin{example}\label{ex:hibi}
	Let $L_0$ be the planar distributive lattice of Figure~\ref{fig:hibi}, and for a square $[x,y]$ of $L_0$, that is, an interval with two atoms, let $L$ be the lattice obtained from $L_0$ by adding a new element $z$ with $x\lessdot z\lessdot y$. Then $L$ is a modular lattice whose only diamond is $[x,y]$, and $z$ is the only element of its fibre. Here $n=11$ and $r=5$, and we consider the four squares of $L_0$.
	\begin{enumerate}
		\item[(1)] The diamond is $[\hat0,c]$, with atoms $a,b,z$. The fibres of the atoms are $\{a\}$, the chain $b<d<f$ and $\{z\}$, so every atom is complemented in its fibre and $I_L$ is radical by Theorem B. However, the closed chain $\{\hat0,a,g,\hat1\}$ belongs to $\cD(L)$ and has rank $3<r-1$. By Theorem A the minimal primes of $I_L$ have heights $6$ and $7$, so $I_L$ is not unmixed.
		\item[(2)] The diamond is $[b,e]$, with atoms $c,d,z$. The fibres are the chains $a<c$ and $d<f$ and $\{z\}$, so $I_L$ is radical. The sets $\{d,e,f,g,h,\hat1\}$ and $\{z,e,g,h,\hat1\}$ belong to $\cD(L)$ and have rank $3$, so again the minimal primes of $I_L$ have heights $6$ and $7$, and $I_L$ is not unmixed.
		\item[(3)] The diamond is $[d,h]$, with atoms $e,f,z$. The fibre of $e$ is the chain $a<c<e<g$, in which $e$ has no complement; this is the configuration of Corollary~\ref{cor:local}. Hence $I_L$ is not radical by Theorem B, and Proposition~\ref{prop:witness} gives $x_{\hat0}x_{\hat1}(x_ex_g-x_dx_{\hat1})\in\sqrt{I_L}\setminus I_L$. Moreover, the closed chains $\{f,h,\hat1\}$ and $\{z,h,\hat1\}$ belong to $\cD(L)$ and have rank $2$, and the minimal primes of $I_L$ have heights $6$, $7$ and $8$.
		\item[(4)] The diamond is $[e,\hat1]$, with atoms $g,h,z$. The fibres are $\{g\}$, the chain $f<h$ and $\{z\}$, so $I_L$ is radical. Moreover, all six members of $\cD(L)$ have rank $4=r-1$, so all minimal primes of $I_L$ have height $6$, and $I_L$ is unmixed by Corollary C.
	\end{enumerate}
\end{example}

Tests with Macaulay2 \cite{M2} on modular lattices, together with the results of this paper, suggest the following conjectures, where $\cD(L)$ is as in Definition~\ref{def:DN}.

\begin{conjecture}\label{conj:unmixed}
For a finite modular lattice $L$ the following are equivalent:
\begin{enumerate}
\item[(a)] $I_L$ is unmixed;
\item[(b)] all minimal primes of $I_L$ have the same height;
\item[(c)] $\rk D=\kappa(L)$ for every $D\in\cD(L)$.
\end{enumerate}
\end{conjecture}

\begin{conjecture}\label{conj:radical}
Let $L$ be a finite modular lattice. If $I_L$ is unmixed, then $I_L$ is radical.
\end{conjecture}

If $L$ has exactly one diamond, then $\kappa(L)=\rk L-1$ by Lemma~\ref{lem:onediamond}, and both conjectures hold by Corollary C. The converse of Conjecture~\ref{conj:radical} is false, as Example~\ref{ex:A} shows. The modular lattices with $\kappa(L)=\rk L-1$ also include lattices with several diamonds, such as $M_3\times\{0<1\}$, and they seem a natural class in which to extend Theorems A and B.

\bigskip

\begin{footnotesize}
	
	{\bf Declaration of competing interest.}  The authors declare that they have no known competing financial interests or personal relationships that could have appeared to influence the work reported in this paper.
	
	\medskip
	
	{\bf Data availability.} No data was used for the research described in the article. 
	
\end{footnotesize}
\bigskip

\end{document}